\input epsf
\magnification=\magstep1
\documentstyle{amsppt}
\pagewidth{6.48truein}\hcorrection{0.0in}
\pageheight{9.0truein}\vcorrection{0.0in}
\TagsOnRight
\NoRunningHeads
\catcode`\@=11
\def\logo@{}
\footline={\ifnum\pageno>1 \hfil\folio\hfil\else\hfil\fi}
\topmatter
\title Rotational invariance of quadromer correlations on the hexagonal lattice 
\endtitle
\author Mihai Ciucu\endauthor
\thanks Research supported in part by NSF grants DMS 9802390 and DMS 0100950.
\endthanks
\affil
  School of Mathematics, Georgia Institute of Technology\\
  Atlanta, Georgia 30332-0160
\endaffil
\abstract
In 1963 Fisher and Stephenson \cite{FS} conjectured that the monomer-monomer
correlation on the square lattice is rotationally invariant. In this paper we
prove a closely related statement on the hexagonal lattice. Namely, we consider 
correlations of two quadromers (four-vertex subgraphs consisting of a monomer and its 
three neighbors) and show that they are rotationally invariant.
\endabstract
\endtopmatter
\document

\def\mysec#1{\bigskip\centerline{\bf #1}\message{ * }\nopagebreak\par\bigskip}

\def\myref#1{\item"{[{\bf #1}]}"} 
 
\def\pf{{\it Proof.\ }} 
\def\endpf{\hbox{\qed}\bigskip}
\def\cite#1{\relaxnext@
  \def\nextiii@##1,##2\end@{[{\bf##1},\,##2]}%
  \in@,{#1}\ifin@\def\next{\nextiii@#1\end@}\else
  \def\next{[{\bf#1}]}\fi\next}
\def\proclaimheadfont@{\smc}

\def\pf{{\it Proof.\ }}

\define\Z{{\Bbb Z}}
\define\Q{{\Bbb Q}}
\define\R{{\Bbb R}}
\define\C{{\Bbb C}}
\define\M{\operatorname{M}}
\define\wt{\operatorname{wt}}

\define\twoline#1#2{\line{\hfill{\smc #1}\hfill{\smc #2}\hfill}}

\def\mypic#1{\epsffile{figs/#1}}



\mysec{1. Introduction}

The lines of the square grid divide the plane into unit squares called {\it monomers}.
A {\it dimer} is the union of two monomers that share an edge, and a collection
of disjoint dimers is said to be a {\it dimer tiling} of the planar region
obtained as their union. We denote by $\M(R)$ the number of dimer tilings of the
region $R$.

Let $S_n$ be the square $[0,2n]\times[0,2n]$. Denote by $S_n(r,s;p,q)$ the region 
obtained from $S_n$ by removing the two monomers whose lower left corners have 
coordinates $(r,s)$ and $(r+p,s+q)$, $0\leq r,s,r+p,s+q\leq 2n-1$. 

The {\it boundary-influenced} correlation of the two removed monomers is defined in 
\cite{FS} as
$$\omega_b(r,s;p,q):=\lim_{n\to\infty}\frac{\M(S_n(r,s;p,q))}{\M(S_n)}.\tag1.1$$

If the monomers of $S_n$ are colored in a chessboard fashion, every dimer consists of
one monomer of each color, and one sees that
$\omega_b(r,s;p,q)=0$ unless the two removed monomers have opposite colors. Therefore 
in the following discussion we assume that the latter condition holds (this amounts to
$p$ and $q$ having opposite parities).

The correlation at the {\it center} is then defined as
$$\omega(p,q):=\lim_{r,s\to\infty}\omega_b(r,s;p,q).\tag1.2$$

The results in \cite{FS} provide explicit exact values for $\omega(p,q)$ for small,
concrete values of $p$ and $q$. In the cases $q=0$ (or $q=1$) and $p=q+1$---which
correspond to a lattice and a lattice diagonal direction, respectively---, 
Fisher and Stephenson \cite{FS}
computed the values of these correlations for the first thirteen
(respectively, first 22)
values of $p$. These tables provide strong evidence that the correlations
$\omega(p,q)$ decay to zero at exactly the same rate (namely, as the inverse
square root of the distance between the monomers, with exactly the same constant
of proportionality) along this two
inequivalent directions\footnote{Hartwig \cite{H} proved that the correlations decay
as predicted by Fisher and Stephenson along the lattice diagonal direction.}. 
After making this remark, Fisher and Stephenson continue in \cite{FS} by stating:
``This equality suggests that the monomer correlations decay
isotropically (with full angular symmetry).'' This conjecture was the starting point
of the present paper.


Clearly, these considerations can also be made on the triangular lattice. A monomer is
then a unit equilateral triangle, and the dimers are the unit rhombi consisting of the 
union of two monomers with a common edge. 

Each dimer covers one up-pointing and one
down-pointing monomer, so the two orientations of monomers play in this case the role
of the two colors in the chessboard coloring of the square lattice.
Therefore, the natural analog of the conjecture of Fisher and Stephenson would be the
statement that the correlation of two monomers of opposite orientation is rotationally
invariant.

A {\it quadromer} is by definition the union
of any monomer with the three monomers adjacent to it. Clearly, this is just a
lattice triangle of side two. 

It is correlations of two quadromers that we prove are rotationally invariant in this
paper. Our method of proof is based on exact counting of dimer coverings of certain
regions on the triangular lattice, and the reduction of the problem to these exact
counts works only when we remove triangles of side two---quadromers---, and not when
monomers are removed. In essence, this is due to the fact that two is an even number,
and this allows to express the relevant quantities as convenient determinants, evaluate
them and obtain the asymptotics of the resulting expressions.

We note that removing a quadromer from a region is equivalent to removing any two of
its outer monomers---such a pair could be called a {\it bimer}. Indeed, any dimer
covering of the region resulting by removal of a bimer must contain the unit rhombus 
fitting in its notch, thus extending it to a quadromer.

In order to define quadromer correlations, we need a family of regions on the
triangular lattice to play the role the squares $S_n$ played on the square lattice. The
choice of these regions is perhaps the single most crucial part in our proof. 

To this end, for any integer $n\geq1$ we define the region $P_n$ as follows. 
Fix a vertex $O$ of the triangular lattice and 
consider the zig-zag lattice path $L_u$ of length $2n$ extending upward from $O$ and 
taking alternate steps northeast and northwest (see Figure~1.1). 
Let $O'$ be the southwestern neighbor 
of $O$, and consider a second zig-zag lattice path $L_d$ of length $2n$, extending 
downward from $O'$ and taking alternate steps southeast and southwest. The union of 
these two lattice paths and the segment $OO'$ forms the eastern boundary of $P_n$. The 
rest of its boundary is defined by moving from the lowest vertex of $L_d$ successively 
$n$ units west, $2n$ units northwest, $2n+1$ units northeast and $n$ units
east, along lattice lines throughout, to arrive at the topmost vertex of $L_u$. 
Our region $P_n$ is the region traced out this way, with the additional requirement 
that the $n$ dimer positions closest to $L_u$ have weight 1/2 (these positions are
indicated by shaded ellipses in Figure 1.1): what this means is that a
dimer covering of $P_n$ using precisely $k$ of these dimers gets weight $1/2^k$, and 
$\M(P_n)$ is the {\it weighted} count of the dimer coverings of $P_n$, i.e., the sum 
of all their weights.

\topinsert
\twoline{\mypic{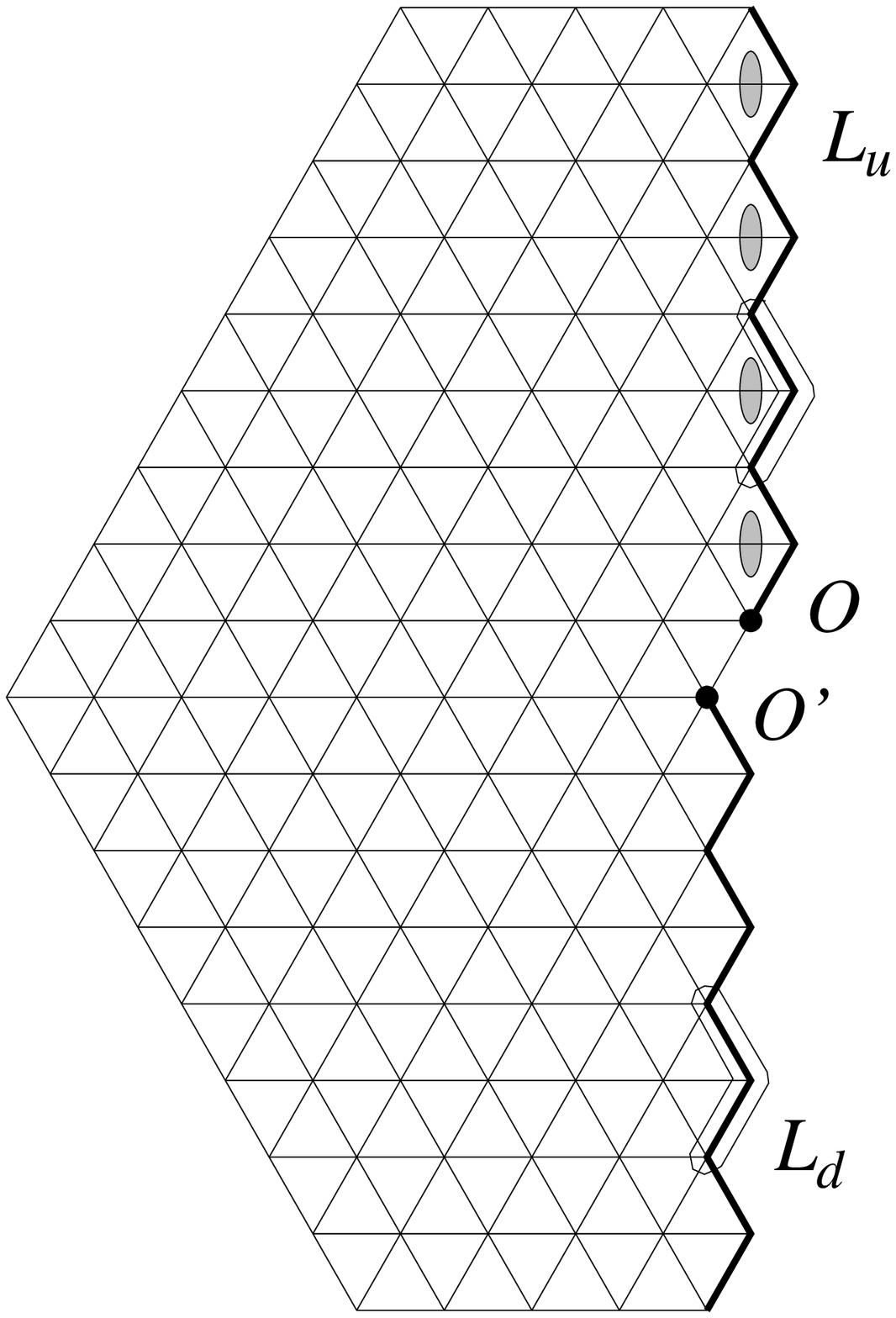}}{\mypic{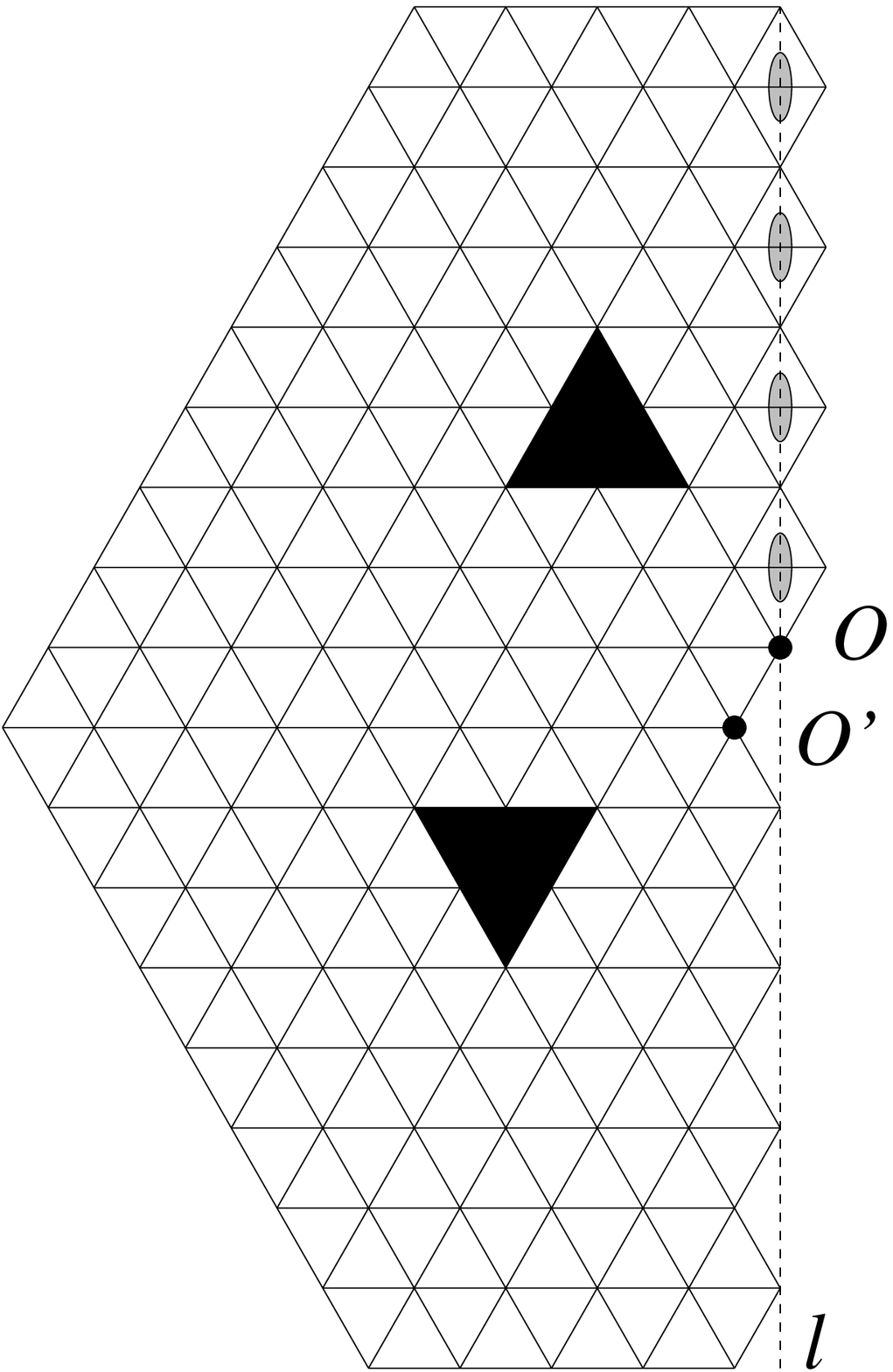}}
\medskip
\twoline{Figure~1.1.{ $P_{4}$.}}
{Figure~1.2.{ $P_{4}(3,0;2,1)$.}}
\endinsert

Let $\ell$ be the vertical line through $O$. 
Denote by $D(R,v)$ the down-pointing quadromer 
whose base is centered $R$ units to the left of $\ell$ and lies on a 
horizontal lattice line crossing $L_d$ $2v+1$ units below $O'$.
Let $U(R,v)$ be the up-pointing quadromer whose base is centered
$R$ units to the left of $\ell$ and lies on a horizontal lattice line crossing $L_u$
$2v$ units above $O$ (Figure 1.2 shows quadromers
$D(3,0)$ and $U(2,1)$)\footnote{Considering the relative
position to $L_u$ and $L_d$, there are two distinct types of both up-pointing
and down-pointing quadromers. The reason we made the indicated choice
out of the total of four possibilities is that it provides even vertical separations
between the bases of the quadromers. The other three choices would entail only minor
alterations to the considerations in this paper.} .

Let $P_n(R_1,v_1;R_2,v_2)$ be the region obtained from $P_n$ by removing the
quadromers $D(R_1,v_1)$ and $U(R_2,v_2)$, where $R_1,R_2\geq1$ and $v_1,v_2\geq0$ (see
Figure 1.2 for an example).

Paralleling (1.1), we define the 
boundary-influenced correlation of two removed qua\-dro\-mers as

$$\omega_b(R_1,v_1;R_2,v_2):=
\lim_{n\to\infty}\frac{\M(P_n(R_1,v_1;R_2,v_2))}{\M(P_n)}.\tag1.3$$

In analogy to (1.2), we define the correlation of the quadromers at the center, with
horizontal separation $r$ and vertical separation $\sqrt{3}u\geq\sqrt{3}$, by
$$\omega(r,u):=
\lim_{R\to\infty}\omega_b(R+r,u-1;R,0).\tag1.4$$

When making the separation $(r,u)$ of the two quadromers
grow to infinity it is natural to do it so that $u=qr+c$, where $q$ and $c$ are fixed
rational numbers. 

The main theorem of this paper is the following.

\proclaim{Theorem 1.1} As $r$ and $u$ approach infinity so that $u=qr+c$, with 
$q\geq0$ and $c$ fixed rational numbers $($$c\geq1$ if $q=0$$)$,
$$\omega(r,u)=\frac{3}{4\pi^2(r^2+3u^2)}+o(r^{-2}).\tag1.5$$
\endproclaim

Since the parenthesis in 
the denominator in (1.5) is just the square of the distance between the centers
of the bases of the removed quadromers, this theorem shows in particular that the 
quadromer correlation $\omega(r,u)$ is rotationally invariant.

\smallpagebreak
\flushpar
{\smc Remark 1.2.} An alternative (and perhaps more natural) way to define the 
correlation of the two removed quadromers would be to consider say regular hexagons $H_n$ of
side $n$ instead of our regions $P_n$. One could even remove the two quadromers from near the center 
of $H_n$, and then the analog of (1.3) would define directly the correlation at the center.

However, when carrying out the line of approach of our proof in this set-up, instead of the 
well-behaved regions $P_n[k_1,k_2;l_1,l_2]$ of Section 2, one is lead to regions obtained from $H_n$
by removing pairs of unit triangles from along two of its opposite sides---and these turn out not to
possess a simple product formula for the number of their tilings. This precludes having an analog of 
Lemma 2.1, which is the first key step in our proof.

\mysec{2. A quadruple sum for quadromer correlations}

The reason we chose the regions $P_n$ as above is that they have, as shown in
\cite{C}, the
remarkable property that quite general alterations of their eastern boundary create
regions whose weighted dimer counts are given by simple product formulas. 

To state this precisely, view the lattice paths $L_u$ and $L_d$ on the eastern
boundary as consisting of $n$ {\it bumps} each---successions of two lattice steps
forming an angle opening to the left (one sample bump on each of $L_d$ and $L_u$ is indicated in
Figure 1.1). Label the bumps on each path successively by
$0,1,\dotsc,n-1$, starting from the bumps closest to $O$\footnote{Note that this
differs from the labeling used in \cite{C}, which is obtained by increasing the
present labels by 1.} (see Figure 2.1). Let $B$ be a bump on $L_d$, and consider the unique
down-pointing quadromer $D$ that contains $B$. By removing bump $B$ from $P_n$ we mean removing the 
three monomers of $D$ contained in $P_n$. We define removal of a bump of $L_u$ analogously, via
the unique {\it up}-pointing quadromer containing it. Figure 2.1 illustrates the effect of removing
bumps 0 and 2 of $L_d$ and bumps 1 and 2 of $L_u$ from the region $P_4$.


\topinsert
\centerline{\mypic{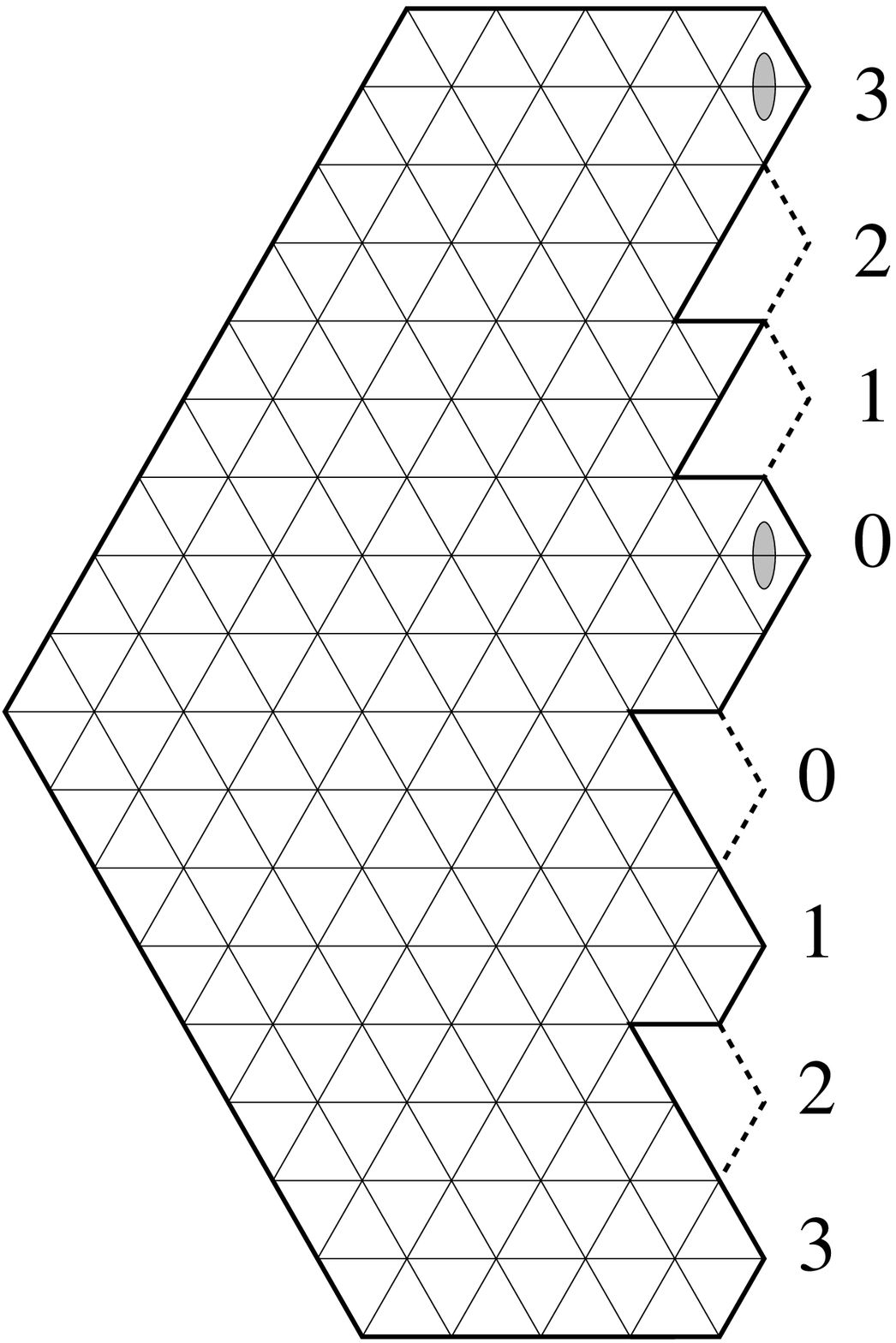}}
\medskip
\centerline{{\smc Figure~2.1.} {\rm $P_4[0,2;1,2]$.}}
\endinsert

Let $P_n[k_1,k_2;l_1,l_2]$ be the region obtained from $P_n$ by removing bumps $k_1$
and $k_2$ from $L_d$ and bumps $l_1$ and $l_2$ from $L_u$, where 
$0\leq k_1<k_2\leq n-1$ and $0\leq l_1<l_2\leq n-1$ (Figure 2.1 shows $P_4[0,2;1,2]$). 
In \cite{C,(2.2),\,(1.1)--(1.6)} explicit
simple product formulas (i.e., with factors of size at most linear in the parameters) 
are given for the weighted count of dimer coverings of a family of regions
$\bar{R}_{\bold l,\bold q}(x)$ that includes the $P_n[k_1,k_2;l_1,l_2]$'s 
(${\bold l}$ and ${\bold q}$ are lists of strictly increasing positive integers, 
and $x$ is a nonnegative integer; in the notation of \cite{C}, 
$\bar{R}_{[1,\dotsc,k_1,k_1+2,\dotsc,k_2,k_2+2,\dotsc,n],
[1,\dotsc,l_1,l_1+2,\dotsc,l_2,l_2+2,\dotsc,n]}(n)$ is the
region that we denote here $P_n[k_1,k_2;l_1,l_2]$). Thus we obtain simple product
formulas for the numbers $\M(P_n[k_1,k_2;l_1,l_2])$. Using them, one obtains after
straightforward if somewhat lengthy manipulations that
$$
\align
\lim_{n\to\infty}\frac{\M(P_n[k_1,k_2;l_1,l_2])}{\M(P_n)}
=&\frac{2^{-4}(2k_1+1)!\,(2k_2+1)!\,(2l_1+1)!\,(2l_2+1)!}
{2^{2k_1+2k_2+2l_1+2l_2}k_1!\,(k_1+1)!\,k_2!\,(k_2+1)!\,l_1!^2\,l_2!^2}\\
\times&\frac{(k_2-k_1)(l_2-l_1)}{(k_1+l_1+2)(k_1+l_2+2)(k_2+l_1+2)(k_2+l_2+2)}.\tag2.1
\endalign
$$
\proclaim{Lemma 2.1} The boundary-influenced correlation $\omega_b(R_1,v_1;R_2,v_2)$
is given by
$$\align
\omega_b(R_1,v_1;R_2,v_2)=&2^{-4}R_1R_2(R_2-1/2)(R_2+1/2)\\
\times|\sum_{a,b=0}^{R_1}\sum_{c,d=0}^{R_2}&(-1)^{a+b+c+d}
\frac{(R_1+a-1)!\,(R_1+b-1)!}
{(2a)!\,(R_1-a)!\,(2b)!\,(R_1-b)!}\\
\times&\frac{(R_2+c-1)!\,(R_2+d-1)!}
{(2c+1)!\,(R_2-c)!\,(2d+1)!\,(R_2-d)!}\\
\times&\frac{(2v_1+2a+1)!\,(2v_1+2b+1)!}
{2^{2(2v_1+a+b)}(v_1+a)!\,(v_1+a+1)!\,(v_1+b)!\,(v_1+b+1)!}\\
\times&\frac{(2v_2+2c+1)!\,(2v_2+2d+1)!}
{2^{2(2v_2+c+d)}(v_2+c)!^2\,(v_2+d)!^2}\\
\times&\frac{(b-a)^2(d-c)^2}
{(u+a+c)(u+a+d)(u+b+c)(u+b+d)}|,\tag2.2
\endalign$$
where $u=v_1+v_2+2$.
\endproclaim

To prove this Lemma we will need the following special case of the 
Lindstr\"om-Gessel-Viennot theorem on non-intersecting lattice paths. 

Our lattice paths will be paths on the directed grid graph $\Z^2$, with edges 
oriented so that they point in the positive direction. We allow the edges of $\Z^2$ to 
be weighted, and define the weight of a lattice path to be the product of the weights 
on its steps. The weight of an $N$-tuple of lattice paths is the product of the
individual weights of its members. 
The weighted count of a set of $N$-tuples of 
lattice paths is the sum of the weights of its elements.

Let ${\bold u}=(u_1,\dotsc,u_N)$ and ${\bold v}=(v_1,\dotsc,v_N)$ be two fixed sets of
starting and ending points on $\Z^2$, and let  $\Cal N({\bold u},{\bold v})$ be the set 
of non-intersecting lattice paths with these starting and ending points. For 
${\bold P}\in\Cal N({\bold u},{\bold v})$, let $\sigma_{\bold P}$ be the permutation 
induced by ${\bold P}$ on the set consisting of the $N$ indices of its starting and 
ending points.

\proclaim{Theorem 2.2 (Lindstr\"om-Gessel-Viennot \cite{GV})} 
$$\sum_{{\bold P}\in\Cal N({\bold u},{\bold v})}(-1)^{\sigma_{\bold P}}\wt({\bold P})=
\det\left((a_{ij})_{1\leq i,j\leq n}\right),$$
where $a_{ij}$ is the weighted count of the lattice paths from $u_i$ to $v_j$.
\endproclaim

What makes possible the use of this result in our setting is a well-known procedure of
encoding dimer coverings by families of
non-intersecting ``paths of dimers:'' given a dimer covering $T$ of a region $R$ on the 
triangular lattice and a lattice line direction $d$, the dimers of $T$ parallel to $d$
(i.e., having two sides parallel to $d$) can naturally be grouped into
non-intersecting paths joining the lattice segments on the boundary of $R$ that are
parallel to $d$, and conversely this family of paths determines the dimer covering
(see Figure 2.2 for an illustration of this and e.g. \cite{C} for a more detailed 
account). 

We will find it convenient to view the paths of dimers directly as lattice paths on 
$\Z^2$, thus bypassing the ``extra steps'' of bijecting them with lattice paths on a 
lattice of rhombi with angles of 60 and 120 degrees, and then deforming this to the 
square lattice. In this context, the ``points'' of our lattice $\Cal L$ are the 
edges of the triangular lattice parallel to a chosen lattice line direction $d$---we 
call them {\it segments}---, and the ``lines''of $\Cal L$  are sequences of adjacent 
dimers extending along the two lattice line directions different from $d$. The dimers 
that these lines consist of are the {\it edges} of $\Cal L$. 

\smallpagebreak
{\it Proof of Lemma 2.1.}  
Choose the lattice line direction $d$ in the above encoding procedure to be the
southwest-northeast direction, and choose the positive directions in the lattice 
$\Cal L$ so that they point
east and southeast. Encode the dimer coverings of
$P_n(R_1,v_1;R_2,v_2)$ by $(2n+3)$-tuples of non-intersecting paths consisting of 
dimers parallel to $d$ (see Figure 2.2; there and in the following figures the dimer
positions weighted by 1/2 are not distinguished, but are understood to carry that
weight). 

Let ${\bold P}$ be such a $(2n+3)$-tuple. Consider the
permutation $\sigma_{\bold P}$ induced by ${\bold P}$ on the set of the $2n+3$ indices 
of its starting and ending points. We claim that the sign of $\sigma_{\bold P}$ is
independent of ${\bold P}$. Indeed, denote by $U$ and $D$ 
the removed up-pointing and down-pointing quadromers, respectively. While the way in
which the starting and ending points of ${\bold P}$---clearly independent of ${\bold
P}$---are matched up depends on ${\bold P}$, it is always the case that the two paths
ending on $U$ start at consecutive starting points, and the two starting at $D$ end at
consecutive ending points. It is easy to see that this implies that all 
$\sigma_{\bold P}$'s have the same sign.

\topinsert
\twoline{\mypic{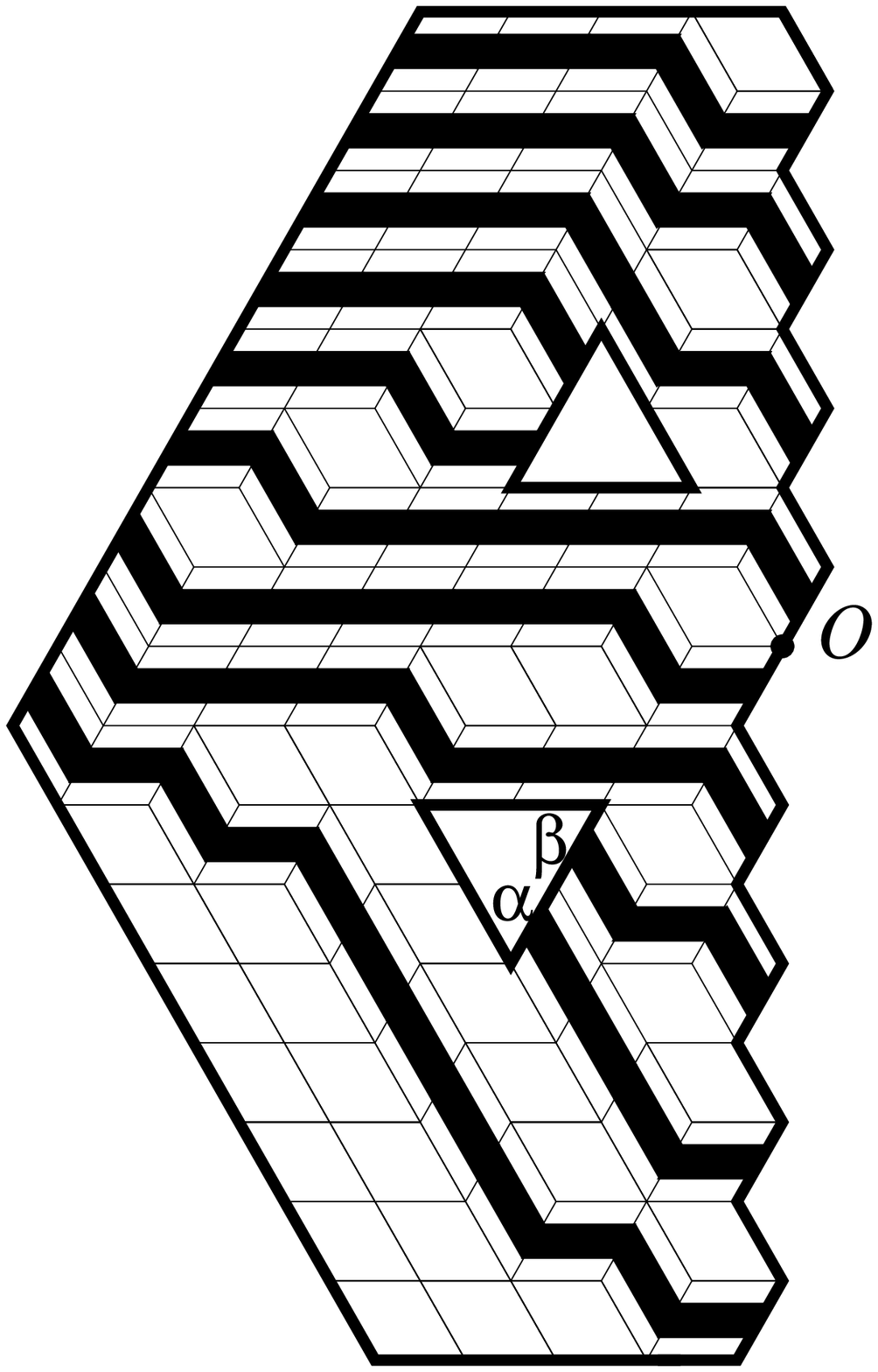}}{\mypic{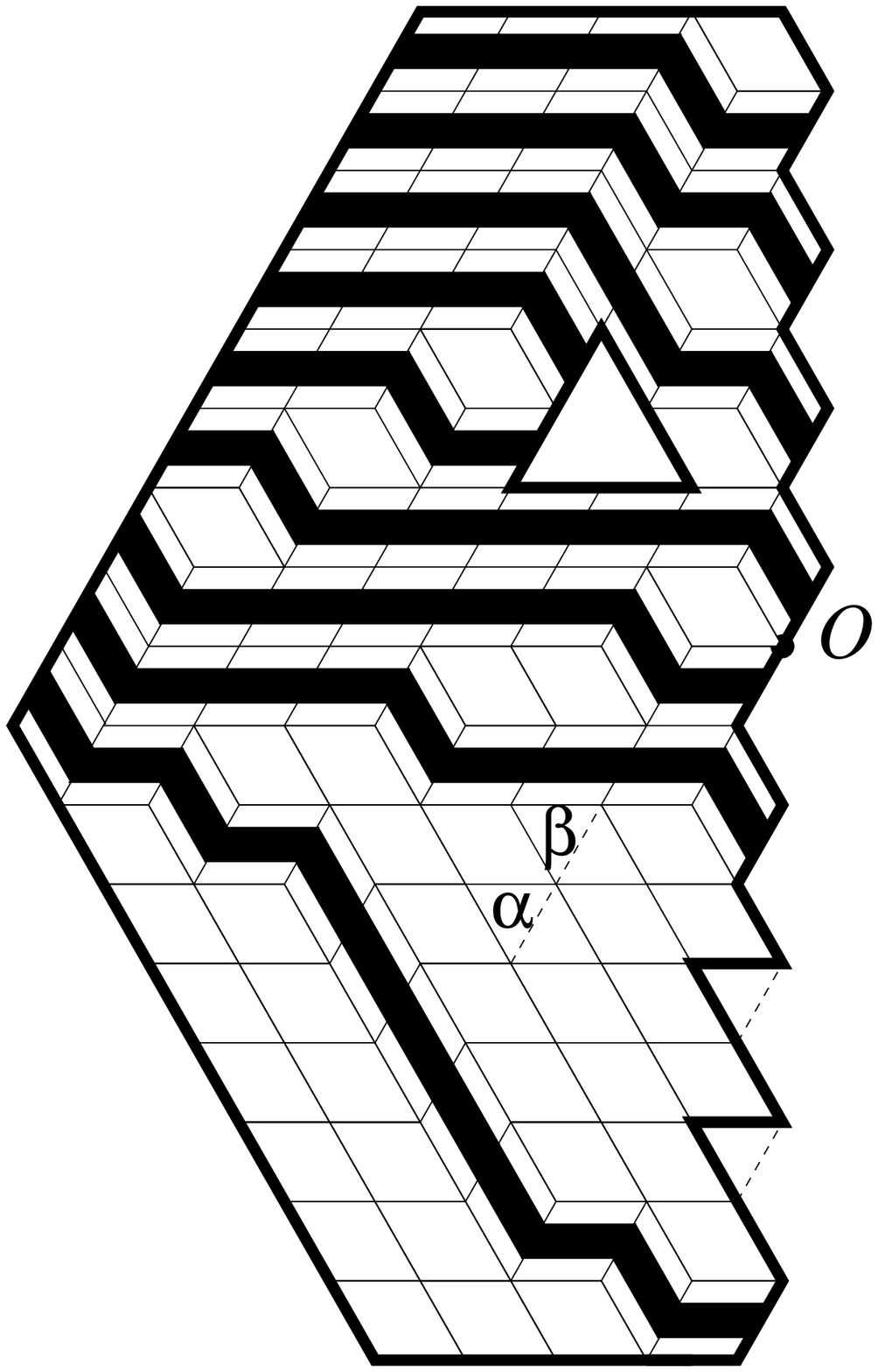}}
\medskip
\twoline{Figure~2.2.{ }}
{Figure~2.3.{ }}
\twoline{{ \rm A dimer covering of $P_4(3,0;2,1)$}}
{{ \rm The effect of Laplace expansion over}}
\twoline{{\rm encoded by paths of rhombi.}}
{{\rm the rows indexed by $\alpha$ and $\beta$.}}
\endinsert

Weight by 1/2 the edges of our ``path-encoding'' lattice $\Cal L$ corresponding to
dimer positions weighted by 1/2 in $P_n(R_1,v_1;R_2,v_2)$.  Weight all other edges of 
$\Cal L$ by 1.
Then the weight of the dimer covering encoded by ${\bold P}$ is just $\wt({\bold P})$,
and we obtain by Theorem 2.2 and the constancy of the sign of $\sigma_{\bold P}$ that
$$\M(P_n(R_1,v_1;R_2,v_2))=\left|\det A\right|,\tag2.3$$
where $A$ is the $(2n+3)\times(2n+3)$ matrix recording the weighted counts of the
lattice paths with given starting and ending points (note that the right hand side of
(2.3) is independent of the ordering of these starting and ending points).

We deduce (2.2) by applying Laplace expansion to the determinant in (2.3). Recall that
for any $m\times m$ matrix $M$ and any $s$-subset $S$ of $[m]:=\{1,\dotsc,m\}$, 
Laplace expansion along the rows with indices in $S$ states that
$$\det M=
\sum_K(-1)^{\epsilon(K)}\det M_S^K\det M_{[m]\setminus S}^{[m]\setminus K},\tag2.4$$
where $K$ ranges over all $s$-subsets of $[m]$, 
$\epsilon(K):=\sum_{k\in K}(k-1)$ and $M_I^J$ is the 
submatrix of $M$ with row-index set $I$ and column-index set $J$.

The rows and columns of the matrix $A$ in (2.3) are indexed by the starting and ending
points of the $(2n+3)$-tuples of non-intersecting lattice paths encoding the dimer
coverings of $P_n(R_1,v_1;R_2,v_2)$. The starting points are the $2n+1$ unit segments
along the northwestern boundary of $P_n(R_1,v_1;R_2,v_2)$, together with the two
segments $\alpha$ and $\beta$ of $D$ parallel to $d$ (see Figure 2.2). 
The ending points are the $2n+1$
segments parallel to $d$ on the eastern boundary of $P_n(R_1,v_1;R_2,v_2)$, together
with two more such segments on $U$. 

Apply Laplace expansion to the matrix $A$ of (2.3) along the two rows indexed by 
$\alpha$ and $\beta$ (see Figure 2.2). 
The first determinant in the summand in (2.4) is then just a two
by two determinant. Its entries are weighted counts of lattice paths on $\Cal L$ that
start at $\alpha$ or $\beta$ and end at some segment on $L_d$. There are only
$R_1+1$ segments on $L_d$ that can be reached this way. Label them consecutively from
top to bottom by $0,1,\dotsc,R_1$. We can restrict summation in (2.4) to the
two-element subsets $K$ of this set of segments: all other terms have at least
one zero column in the two by two determinant. Therefore we obtain from (2.3) that
$$\M(P_n(R_1,v_1;R_2,v_2))=\left|\sum_{0\leq a,b\leq R_1}(-1)^{a+b}
\det A_{\{\alpha,\beta\}}^{\{a,b\}}
\det A_{[2n+3]\setminus\{\alpha,\beta\}}^{[2n+3]\setminus\{a,b\}}\right|.\tag2.5$$
Choosing the origin of $\Cal L$ to be at $\alpha$, one sees that $\beta$ has 
coordinates $(-1,1)$
and the segment labeled $j$ on $L_d$ has coordinates $(R_1-1-j,2j)$, $j=0,\dotsc,R_1$.
Since the lattice paths counted by the entries of $A_{\{\alpha,\beta\}}^{\{a,b\}}$
have all steps weighted by 1, the determinant of this matrix is
$$
\det A_{\{\alpha,\beta\}}^{\{a,b\}}=\det\left[\matrix
{R_1-1+a\choose 2a} {R_1-1+b\choose 2b}\\
{R_1-1+a\choose 2a-1} {R_1-1+b\choose 2b-1}
\endmatrix\right]=
2R_1\frac{(b-a)(R_1+a-1)!\,(R_1+b-1)!}{(2a)!\,(R_1-a)!\,(2b)!\,(R_1-b)!}.\tag2.6
$$
On the other hand, the second determinant in the summand in (2.5) can be interpreted
as being the weighted count of dimer coverings of the region
$P_n^{[a,b]}(R_2,v_2)$ obtained from
$P_n(R_1,v_1;R_2,v_2)$ by placing back quadromer D and removing the two monomers (i.e., unit 
triangles) that
contain the segments labeled $a$ and $b$ on $L_d$, under the labeling of the preceding paragraph 
(see Figure 2.3 for an illustration). Indeed, the
Lindstr\"om-Gessel-Viennot matrix of this region is precisely
$A_{[2n+3]\setminus\{\alpha,\beta\}}^{[2n+3]\setminus\{a,b\}}$, and by the argument 
that proved (2.3) we
obtain that $\M(P_n^{[a,b]}(R_2,v_2))$ is equal to
$\det A_{[2n+3]\setminus\{\alpha,\beta\}}^{[2n+3]\setminus\{a,b\}}$, up to a sign that 
is independent of $a$ and $b$ (indeed, the permutations $\sigma_{\bold P}$ that occur
when applying Theorem 2.2 to the region $P_n^{[a,b]}(R_2,v_2)$ are independent of $a$
and $b$).
Therefore, using (2.6) we can rewrite (2.5) as
$$\align
\M(P_n&(R_1,v_1;R_2,v_2))=\\
&2R_1\left|\sum_{0\leq a<b\leq R_1}(-1)^{a+b}
\frac{(b-a)(R_1+a-1)!\,(R_1+b-1)!}{(2a)!\,(R_1-a)!\,(2b)!\,(R_1-b)!}
\M(P_n^{[a,b]}(R_2,v_2))\right|.\tag2.7
\endalign$$
In turn, $\M(P_n^{[a,b]}(R_2,v_2))$ can be expressed by a formula similar to the one 
above. To obtain this, encode the tilings of $P_n^{[a,b]}(R_2,v_2)$ by lattice paths,
choosing this time the lattice direction $d$ to be the southeast-northwest direction,
and the positive directions in the encoding lattice $\Cal L$ to point east and
northeast. As in the previous encoding, weight by 1/2 those segments of $\Cal L$ that
correspond to dimer positions weighted 1/2 in $P_n^{[a,b]}(R_2,v_2)$, and weight all
its remaining segments by 1.

Each tiling of $P_n^{[a,b]}(R_2,v_2)$ gets encoded this way by a
$(2n+2)$-tuple ${\bold P}$ of non-intersecting lattice paths on $\Cal L$, starting at
the unit segments on its southwestern boundary or at the unit segments $\gamma$
and $\delta$ of $U$ that are parallel to $d$, and ending at the unit segments
parallel to $d$ on its eastern boundary (see Figure 2.4).  
As in the argument that proved (2.3), the sign of the permutation $\sigma_{\bold P}$
is independent of ${\bold P}$. Therefore, we obtain by Theorem 2.2 that
$$\M(P_n^{[a,b]}(R_2,v_2))=\epsilon\det B,\tag2.8$$
where $B$ is the $(2n+2)\times(2n+2)$ matrix recording the weighted counts of the
lattice paths with specified starting and ending points, and the sign $\epsilon$ in 
front of the determinant is the same for all choices of $a$ and $b$.

\topinsert
\twoline{\mypic{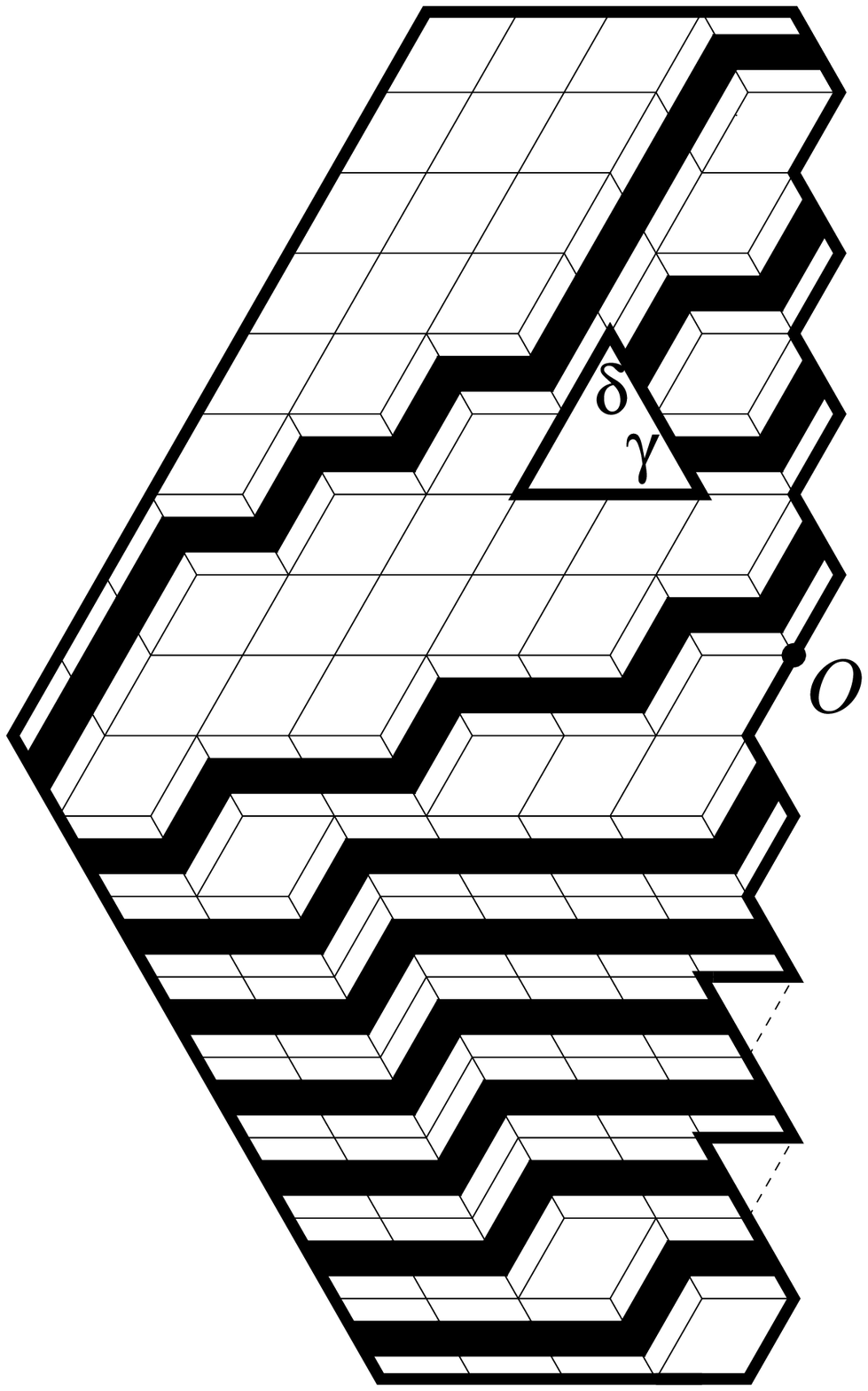}}{\mypic{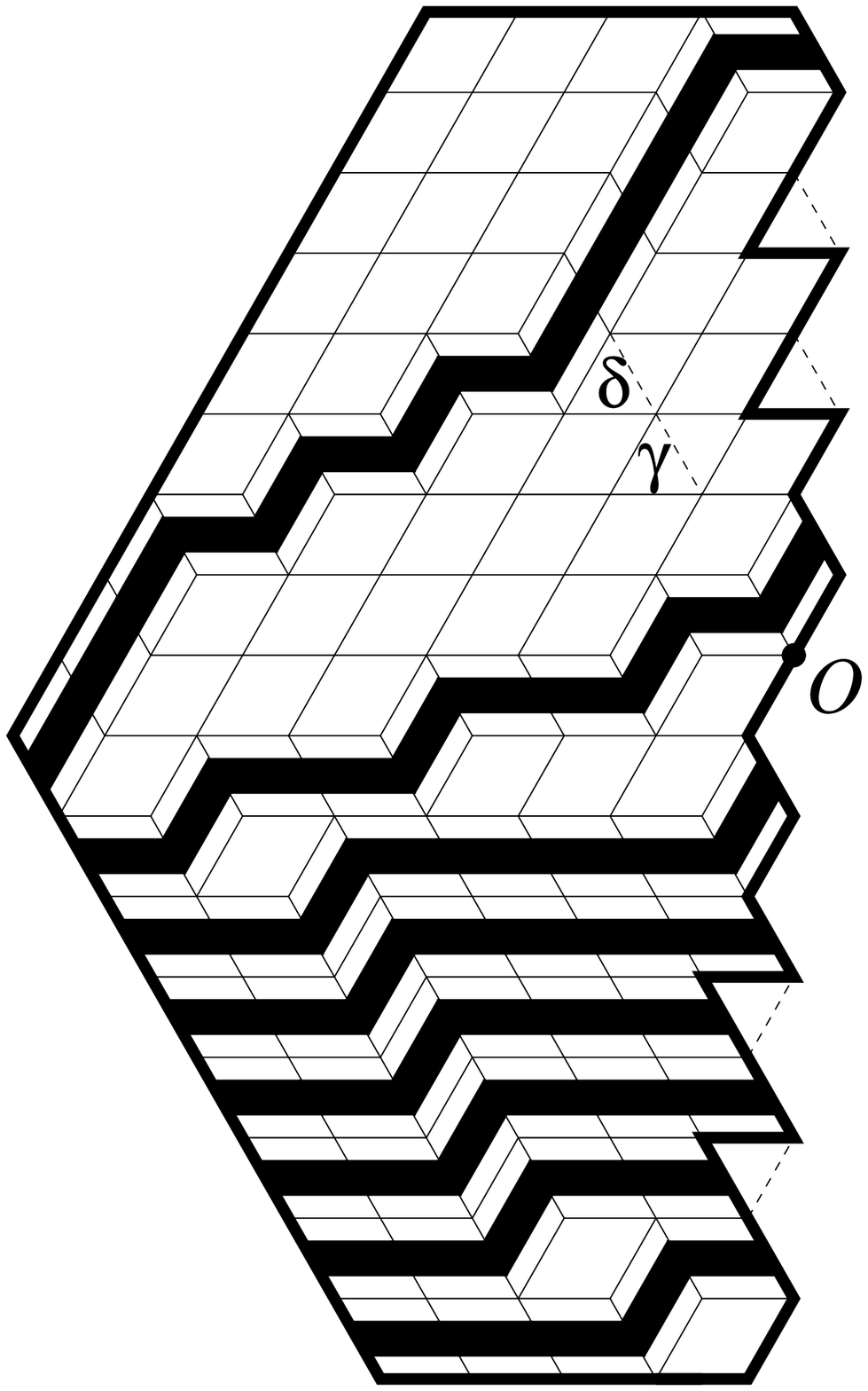}}
\medskip
\twoline{Figure~2.4.{ }}
{Figure~2.5.{ }}
\twoline{{ \rm A dimer covering of $P_4^{[1,2]}(2,1)$}}
{{ \rm The effect of Laplace expansion over}}
\twoline{{\rm \!\!\!\!\! encoded by paths of rhombi.}}
{{\rm the rows indexed by $\gamma$ and $\delta$.}}

\endinsert

Apply Laplace expansion in $\det B$ along the two rows indexed by $\gamma$ and
$\delta$. The first determinant in the summand of (2.4) is again two by two, and
records weighted counts of lattice paths starting at $\gamma$ or $\delta$ and ending
at some segment on $L_u$ (see Figure 2.4). 
There are $R_2+1$ segments on $L_u$ that can be
reached this way; label them consecutively from bottom to top by $0,1,\dotsc,R_2$. 
As with our previous Laplace expansion, we can restrict the summation range in (2.4)
to obtain
$$\M(P_n^{[a,b]}(R_2,v_2))=\epsilon\sum_{0\leq c<d\leq R_1}(-1)^{c+d-1}
\det B_{\{\gamma,\delta\}}^{\{c,d\}}
\det B_{[2n+2]\setminus\{\gamma,\delta\}}^{[2n+2]\setminus\{c,d\}}.\tag2.9$$

Centering $\Cal L$ at $\gamma$, $\delta$ has coordinates $(-1,1)$ and the segment
labeled $j$ on $L_u$ has coordinates $(R_2-1-j,2j+1)$, $j=0,1,\dotsc,R_2$. The
weighted counts involved in the entries of $B_{\{\gamma,\delta\}}^{\{c,d\}}$ 
(which involve this time some steps of weight 1/2) are easily calculated and one obtains
$$\align
\det B_{\{\gamma,\delta\}}^{\{c,d\}}=&\det\left[\matrix
\frac{1}{2}{R_2-1+c\choose 2c}+{R_2-1+c\choose 2c+1} 
\ \ \frac{1}{2}{R_2-1+d\choose 2d}+{R_2-1+d\choose 2d+1}\\
\frac{1}{2}{R_2-1+c\choose 2c-1}+{R_2-1+c\choose 2c} 
\ \ \frac{1}{2}{R_2-1+d\choose 2d-1}+{R_2-1+d\choose 2d}
\endmatrix\right]\\
=&2R_2(R_2-1/2)(R_2+1/2)\frac{(d-c)(R_2+c-1)!\,(R_2+d-1)!}
{(2c+1)!\,(R_2-c)!\,(2d+1)!\,(R_2-d)!}.\tag2.10
\endalign$$
On the other hand, by applying Theorem 2.2 one more time one sees that 
$$\det B_{[2n+2]\setminus\{\gamma,\delta\}}^{[2n+2]\setminus\{c,d\}}=
\epsilon'\M(P_n^{[a,b][c,d]}),\tag2.11$$ 
where $P_n^{[a,b][c,d]}$ is the region obtained from $P_n^{[a,b]}(R_2,v_2)$ by placing
back quadromer $U$ and removing the two monomers near $L_u$ containing segments $c$
and $d$, and the sign $\epsilon'$ is independent of $c$ and $d$. Furthermore, 
$P_n^{[a,b][c,d]}$ differs from the region
$P_n[v_1+a,v_1+b;v_2+c,v_2+d]$ considered at the beginning of this section
only in that the former contains four more dimers, which are weighted by 1 and forced 
to be part of all its dimer coverings (see Figure 2.6); so the two regions have equal 
weighted counts of dimer coverings. Therefore, by (2.7), (2.9), (2.10) and (2.11) we
obtain that

\topinsert
\centerline{\mypic{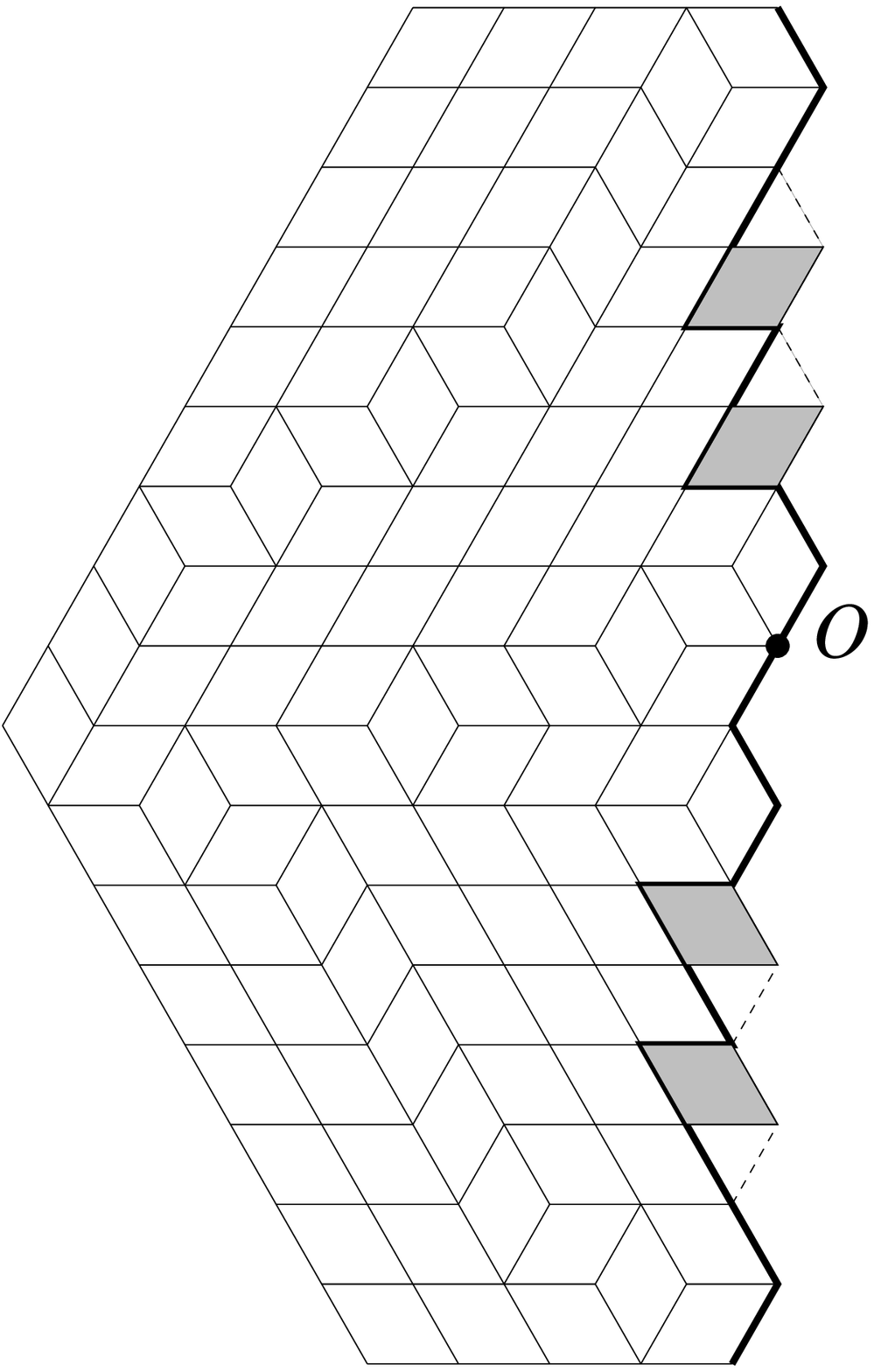}}
\medskip
\centerline{{\smc Figure~2.6.}} 
\centerline{{ \rm The regions $P_4^{[1,2][0,1]}$ and $P_4[1,2;1,2]$}}
\centerline{{ \rm differ only by four dimers.}}
\endinsert
$$
\align
\M(P_n&(R_1,v_1;R_2,v_2))=4R_1R_2(R_2-1/2)(R_2+1/2)\\
&\times|\sum_{0\leq a<b\leq R_1}\sum_{0\leq c<d\leq R_1}
(-1)^{a+b+c+d}
\frac{(b-a)(R_1+a-1)!\,(R_1+b-1)!}{(2a)!\,(R_1-a)!\,(2b)!\,(R_1-b)!}\\
&\ \ \ \ \ \ \ \ \ \ \ \ \ \times
\frac{(d-c)(R_2+c-1)!\,(R_2+d-1)!}{(2c+1)!\,(R_2-c)!\,(2d+1)!\,(R_2-d)!}\\
&\ \ \ \ \ \ \ \ \ \ \ \ \ \times\M(P_n[v_1+a,v_1+b;v_2+c,v_2+d])|.\tag2.12
\endalign
$$
Dividing (2.12) by $\M(P_n)$, letting $n\to\infty$ and using (2.1), one
obtains an expression for $\omega_b(R_1,v_1;R_2,v_2)$ as a quadruple sum in which the
summation indices need to satisfy $a<b$ and $c<d$. 

The fortunate situation is that, on
the one hand, when $a=b$ or $c=d$ the summand becomes zero, and on the other, (2.12)
and (2.1) combine to produce a summand which is invariant under
independently transposing $a$ with $b$ and $c$ with $d$ (because the differences $b-a$
and $d-c$ end up appearing at the second power). Therefore the summation range may
be extended to the one shown in (2.2), at the expense of a multiplicative factor of
$1/4$. 
This leads precisely to the quadruple sum given in the statement of the Lemma. \endpf

\mysec{3. Four double sums, with integral representations}

A simple partial fraction decomposition of part of the summand in (2.2) affords a
great deal of simplification in the expression (2.2) of the boundary-influenced
correlation.

Indeed, one readily checks that
$$\align
&\frac{(b-a)(d-c)}{(u+a+c)(u+a+d)(u+b+c)(u+b+d)}\\
&\ \ \ \ \ \ \ \ \ \ \ \ \ \ =\frac{1}{(u+a+c)(u+b+d)}-\frac{1}{(u+a+d)(u+b+c)}.
\endalign$$
Using this, the portion of the summand contained in the last line of (2.2) becomes
$$\align
&\frac{(b-a)(d-c)}{(u+a+c)(u+b+d)}-\frac{(b-a)(d-c)}{(u+a+d)(u+b+c)}\\
&\ \ \ \ \ \ \ \ \ \ \ \ \ \  
=\frac{ac+bd-ad-bc}{(u+a+c)(u+b+d)}-\frac{ac+bd-ad-bc}{(u+a+d)(u+b+c)}.\tag3.1
\endalign$$
Using this, the fourfold sum of (2.2) becomes a difference of two fourfold sums. The
advantage of this expression is that each of the latter two fourfold sums
can be written as the product of two double sums. Indeed, in all the factors of the
summand in (2.2) except the last line of (2.2), the summation variables can be
separated. Furthermore, the variables $\{a,c\}$ can be separated from $\{b,d\}$
in the first term on the right hand side of (3.1), while $\{a,d\}$ can be separated
from $\{b,c\}$ in the second term of the right hand side of (3.1). The double sums
arising this way are all of the form
$$\align
M_\nu(R_1,R_2):=
\sum_{a=0}^{R_1}\sum_{c=0}^{R_2}&(-1)^{a+c}\frac{(R_1+a-1)!}{(2a)!\,(R_1-a)!}
\frac{(R_2+c-1)!}{(2c+1)!\,(R_2-c)!}\\
\times&\frac{2v_1+2a+1)!}{2^{2v_1+2a}(v_1+a)!\,(v_1+a+1)!}
\frac{2v_2+2c+1)!}{2^{2v_2+2c}(v_2+c)!^2}\frac{\nu}{u+a+c},\tag3.2
\endalign$$
where $\nu$ has one of the values 1, $a$, $c$ and $ac$. For notational convenience
we will often write simply $M_\nu$ instead of $M_\nu(R_1,R_2)$.

More precisely, consider the term $ac$ of the numerator of the first fraction on the
right hand side of (3.1). When summing over $a$, $b$, $c$ and $d$ as required by
(2.2), this term gives rise to $M_{ac}M_1$. Similarly, the remaining terms of that
numerator, $bd$, $-ad$ and $-bc$, give rise to $M_1M_{ac}$, $-M_aM_c$ and $-M_aM_c$,
respectively. In the same fashion, the second term on the right hand side of (3.1)
generates the products $-M_aM_c$, $-M_aM_c$, $M_{ac}M_1$ and $M_1M_{ac}$, respectively.
Therefore, we obtain by (2.2) that
$$
\align
\omega_b(R_1,v_1;R_2,v_2)&=2^{-4}R_1R_2(R_2-1/2)(R_2+1/2)|4M_1M_{ac}-4M_aM_c|\\
&=2^{-2}R_1R_2(R_2-1/2)(R_2+1/2)|M_1M_{ac}-M_aM_c|.
\tag3.3
\endalign
$$
Thus, the asymptotic study of the correlation reduces to studying the asymptotics
of these four double sums. By their definition (3.2), had it not been for the factor
$1/(u+a+c)$, these double sums would further be separable as products of simple sums.
We can get around the obstacle posed by this factor by expressing it as an
integral
\footnote{This useful trick was pointed out to the author independently by
Ira Gessel and Doron Zeilberger. Its use allows a shorter proof for the asymptotics
of the four double sums than our original proof, which relied on the expansion
$$\frac{1}{u+a+b}=\frac{1}{u}\left\{1-\frac{a+c}{u}+\frac{(a+c)(a+c-1)}{u(u+1)}-
\frac{(a+c)(a+c-1)(a+c-2)}{u(u+1)(u+2)}+\cdots\right\}.$$}:
$$\frac{1}{u+a+c}=\int_0^1x^{u+a+c-1}dx.\tag3.4$$
Indeed, substituting this into (3.2), we obtain that for instance $M_1$ is expressed as
$$\align
M_1=\int_0^1&\left(\sum_{a=0}^{R_1}(-1)^a\frac{(R_1+a-1)!}{(2a)!\,(R_1-a)!}
\frac{(2v_1+2a+1)!}{2^{2v_1+2a}(v_1+a)!\,(v_1+a+1)!}x^a\right)\\
\times&\left(\sum_{c=0}^{R_2}(-1)^c\frac{(R_2+c-1)!}{(2c+1)!\,(R_2-c)!}
\frac{(2v_2+2c+1)!}{2^{2v_2+2c}(v_2+c)!^2}x^c\right)x^{u-1}dx.\tag3.5
\endalign$$
It is not difficult (indeed, with access to a computer algebra package like Maple, it
is immediate) to see that the two sums in the above integral can be expressed in
terms of hypergeometric functions\footnote{The hypergeometric function of parameters
$a_1,\dotsc,a_p$ and $b_1,\dotsc,b_q$ is defined by
$${}_p F_q\!\left[\matrix a_1,\dotsc,a_p\\ b_1,\dotsc,b_q\endmatrix;
z\right]=\sum _{k=0} ^{\infty}\frac {(a_1)_k\cdots(a_p)_k}
{k!\,(b_1)_k\cdots(b_q)_k} z^k\ ,$$
where $(a)_0:=1$ and $(a)_k:=a(a+1)\cdots (a+k-1)$ for $k\geq1$.
} as
$$\align
\sum_{a=0}^{R_1}(-1)^a\frac{(R_1+a-1)!}{(2a)!\,(R_1-a)!}
&\frac{(2v_1+2a+1)!}{2^{2v_1+2a}(v_1+a)!\,(v_1+a+1)!}x^a\\
&=\frac{1}{R_1}\frac{(2v_1+1)!}{2^{2v_1}v_1!\,(v_1+1)!}
\,{} _3 F_2\!\left[\matrix{-R_1,R_1,v_1+{3\over2}}\\{{1\over2},v_1+2}\endmatrix; 
{\displaystyle {x\over4}}\right]\tag3.6
\endalign$$
and
$$\align
\sum_{c=0}^{R_2}(-1)^c\frac{(R_2+c-1)!}{(2c+1)!\,(R_2-c)!}
&\frac{(2v_2+2c+1)!}{2^{2v_2+2c}(v_1+c)!^2}x^c\\
&=\frac{1}{R_2}\frac{(2v_2+1)!}{2^{2v_2}v_2!^2}
\,{} _3 F_2\!\left[\matrix{-R_2,R_2,v_2+{3\over2}}\\{{3\over2},v_2+1}\endmatrix; 
{\displaystyle {x\over4}}\right].\tag3.7
\endalign$$
Replacing these formulas in (3.5) we obtain the following result.

\proclaim{Proposition 3.1} The double sum $M_1$ has the integral representation
$$\align
M_1=\frac{1}{R_1R_2}&\frac{(2v_1+1)!\,(2v_2+1)!}{2^{2v_1+2v_2}v_1!\,(v_1+1)!\,v_2!^2}\\
\times&\int_0^1 
{} _3 F_2\!\left[\matrix{-R_1,R_1,v_1+{3\over2}}\\{{1\over2},v_1+2}\endmatrix; 
{\displaystyle {x\over4}}\right]
\,{} _3 F_2\!\left[\matrix{-R_2,R_2,v_2+{3\over2}}\\{{3\over2},v_2+1}\endmatrix; 
{\displaystyle {x\over4}}\right]x^{u-1}dx.\tag3.8
\endalign
$$
\endproclaim

\smallpagebreak
When applying the same reasoning to the remaining double sums $M_a$, $M_c$ and
$M_{ac}$, two more sums besides (3.6) and (3.7) need to be expressed in terms of
hypergeometric functions. These are
$$\align
\sum_{a=0}^{R_1}(-1)^a&\frac{(R_1+a-1)!}{(2a)!\,(R_1-a)!}
\frac{2v_1+2a+1)!}{2^{2v_1+2a}(v_1+a)!\,(v_1+a+1)!}x^aa\\
&=-R_1x\frac{(2v_1+3)(2v_1+1)!}{2^{2v_1+2}v_1!\,(v_1+2)!}
\,{} _3 F_2\!\left[\matrix{-R_1+1,R_1+1,v_1+{5\over2}}\\{{3\over2},v_1+3}\endmatrix; 
{\displaystyle {x\over4}}\right]\tag3.9
\endalign$$
and
$$\align
\sum_{c=0}^{R_2}(-1)^c&\frac{(R_2+c-1)!}{(2c+1)!\,(R_2-c)!}
\frac{2v_2+2c+1)!}{2^{2v_2+2c}(v_1+c)!^2}x^cc\\
&=-\frac{R_2x}{3}\frac{(2v_2+3)(2v_2+1)!}{2^{2v_2+2}v_2!\,(v_2+1)!}
\,{} _3 F_2\!\left[\matrix{-R_2+1,R_2+1,v_2+{5\over2}}\\{{5\over2},v_2+2}\endmatrix; 
{\displaystyle {x\over4}}\right].\tag3.10
\endalign$$
Using (3.2), (3.4), (3.6), (3.7), (3.9) and (3.10), we obtain the following result.

\proclaim{Proposition 3.2} The double sums $M_a$, $M_c$ and $M_{ac}$ have the integral
representations

$$\align
M_a=-\frac{R_1}{R_2}&\frac{(2v_1+3)(2v_1+1)!\,(2v_2+1)!}
{2^{2v_1+2v_2+2}v_1!\,(v_1+2)!\,v_2!^2}\\
\times&\int_0^1 
{} _3 F_2\!\left[\matrix{-R_1+1,R_1+1,v_1+{5\over2}}\\{{3\over2},v_1+3}\endmatrix; 
{\displaystyle {x\over4}}\right]
\,{} _3 F_2\!\left[\matrix{-R_2,R_2,v_2+{3\over2}}\\{{3\over2},v_2+1}\endmatrix; 
{\displaystyle {x\over4}}\right]x^{u}dx\tag3.11\\ \\
M_c=-\frac{R_2}{3R_1}&\frac{(2v_1+1)!\,(2v_2+1)!\,(2v_2+3)}
{2^{2v_1+2v_2+2}v_1!\,(v_1+1)!\,v_2!\,(v_2+1)!}\\
\times&\int_0^1 
{} _3 F_2\!\left[\matrix{-R_1,R_1,v_1+{3\over2}}\\{{1\over2},v_1+2}\endmatrix; 
{\displaystyle {x\over4}}\right]
\,{} _3 F_2\!\left[\matrix{-R_2+1,R_2+1,v_2+{5\over2}}\\{{5\over2},v_2+2}\endmatrix; 
{\displaystyle {x\over4}}\right]x^{u}dx\tag3.12\\ \\
M_{ac}=\frac{R_1R_2}{3}&\frac{(2v_1+3)(2v_1+1)!\,(2v_2+3)(2v_2+1)!}
{2^{2v_1+2v_2+4}v_1!\,(v_1+2)!\,v_2!\,(v_2+1)!}\\ 
\times\int_0^1& 
{} _3 F_2\!\left[\matrix{-R_1+1,R_1+1,v_1+{5\over2}}\\{{3\over2},v_1+3}\endmatrix; 
{\displaystyle {x\over4}}\right]
\,{} _3 F_2\!\left[\matrix{-R_2+1,R_2+1,v_2+{5\over2}}\\{{5\over2},v_2+2}\endmatrix; 
{\displaystyle {x\over4}}\right]x^{u+1}dx.\tag3.13
\endalign
$$
\endproclaim

\mysec{4. Jacobi polynomials}

The $_3F_2$'s of the preceding section can be expressed in terms of $_2F_1$'s by 
the formula 
$$_3F_2\!\left[\matrix{a,\,b,\,c}\\{a-n,d}\endmatrix;z\right]=\frac{1}{(1-a)_n}
\sum_{k=0}^n (-1)^k{n\choose k}(1-a)_{n-k}\frac{(b)_k(c)_k}{(d)_k}z^k\, 
_2F_1\!\left[\matrix{b+k,c+k}\\{d+k}\endmatrix;z\right],\tag4.1
$$
where $n$ is a nonnegative integer. This follows for instance by comparing the coefficients of the
powers of $z$ and applying the Chu-Vandermonde summation (for the latter, see e.g. 
\cite{GR,(1.2.9),\,p.2}).

Applying this to the $_3F_2$'s in (3.6) and (3.7),
one obtains
$$\align
\!\!\!\!\!\!\!\!\!
&_3F_2\!\left[\matrix{-R_1,R_1,v_1+{3\over2}}\\{{1\over2},v_1+2}\endmatrix;
{\displaystyle {x\over4}}\right]=\frac{1}{(-v_1-{1\over2})_{v_1+1}}\\
\!\!\!\!&\times
\sum_{k=0}^{v_1+1} (-1)^k{v_1+1\choose k}(-v_1-{1\over2})_{v_1+1-k}
\frac{(-R_1)_k(R_1)_k}{(v_1+2)_k}\frac{x^k}{4^k}\, 
_2F_1\!\left[\matrix{-R_1+k,R_1+k}\\{v_1+2+k}\endmatrix;{\displaystyle {x\over4}}
\right]\tag4.2
\endalign
$$
and
$$\align
_3F_2\!&\left[\matrix{-R_2,R_2,v_2+{3\over2}}\\{{3\over2},v_2+1}\endmatrix;
{\displaystyle {x\over4}}\right]=\frac{1}{(-v_2-{1\over2})_{v_2}}\\
&\times
\sum_{l=0}^{v_2} (-1)^l{v_2\choose l}(-v_2-{1\over2})_{v_2-l}
\frac{(-R_2)_l(R_2)_l}{(v_2+1)_l}\frac{x^l}{4^l}\, 
_2F_1\!\left[\matrix{-R_2+l,R_2+l}\\{v_2+1+l}\endmatrix;{\displaystyle {x\over4}}
\right].\tag4.3
\endalign
$$
In turn, the resulting $_2F_1$'s can be expressed in terms of the Jacobi polynomials 
$P_n^{(\alpha,\beta)}(x)$ using
the formula (found for instance in \cite{Sz2,(4.21.2),\,p.62})
$$
_2F_1\!\left[\matrix{n+\alpha+\beta+1,-n}\\{1+\alpha}\endmatrix;
{\displaystyle {1-x\over2}}\right]=\frac{n!\,\Gamma(1+\alpha)}{\Gamma(n+1+\alpha)}
P_n^{(\alpha,\beta)}(x).\tag4.4
$$
We obtain for the $_2F_1$'s of (4.2) and (4.3) the expressions
$$
_2F_1\!\left[\matrix{-R_1+k,R_1+k}\\{v_1+2+k}\endmatrix;
{\displaystyle {x\over4}}\right]=\frac{(R_1-k)!\,(v_1+k+1)!}{(R_1+v_1+1)!}
P_{R_1-k}^{(v_1+k+1,k-v_1-2)}\left(1-{x\over2}\right)\tag4.5
$$
and
$$
_2F_1\!\left[\matrix{-R_2+l,R_2+l}\\{v_2+1+l}\endmatrix;
{\displaystyle {x\over4}}\right]=\frac{(R_2-l)!\,(v_2+l)!}{(R_2+v_2)!}
P_{R_2-l}^{(v_2+l,l-v_2-1)}\left(1-{x\over2}\right).\tag4.6
$$
Substituting (4.2), (4.3), (4.5) and (4.6) in the integral representation (3.8) of
$M_1$, we obtain that
$$
M_1=\sum_{k=0}^{v_1+1}\sum_{l=0}^{v_2} \frac{c_{kl}}{R_1R_2}
\frac{(-R_1)_k(R_1)_k}{(R_1-k+1)_{v_1+k+1}}
\frac{(-R_2)_l(R_2)_l}{(R_2-l+1)_{v_2+l}} I_{kl}(R_1,R_2),\tag4.7
$$
where $c_{kl}$ depends only on $k$, $l$, $v_1$ and $v_2$, and for $0\leq k\leq v_1+1$ 
and $0\leq l\leq v_2$,
$$
I_{kl}(R_1,R_2):=\int_0^1 P_{R_1-k}^{(v_1+k+1,k-v_1-2)}\left(1-{x\over2}\right)
P_{R_2-l}^{(v_2+l,l-v_2-1)}\left(1-{x\over2}\right)x^{k+l+u-1}dx.\tag4.8
$$
By (1.4) and (3.3), we need the asymptotics of $M_1$ for $v_1=u-1$, $v_2=0$, $R_1=R+r$
and $R_2=R$, where $r\geq0$ and $u\geq1$ are fixed and 
$R\to\infty$.

By (4.7), the study of this asymptotics of $M_1$ reduces to the asymptotics of the
$I_{kl}$'s. The asymptotics of the Jacobi polynomials $P_n^{(\alpha,\beta)}$ for large
$n$ is given by the Darboux formula (see e.g. \cite{Sz2, Theorem 8.21.8,\,p.196})
$$\align
P_n^{(\alpha,\beta)}(\cos\theta)&=D_n^{(\alpha,\beta)}(\cos\theta)
+E_n^{(\alpha,\beta)}(\cos\theta)\\
D_n^{(\alpha,\beta)}(\cos\theta)&=
\frac{\cos\left\{\left[n+{\displaystyle{\alpha+\beta+1\over2}}\right]\theta-
\left({\displaystyle{\alpha\over2}+{1\over4}}\right)\pi\right\}}
{\sqrt{\pi n}\left(\sin{\displaystyle{\theta\over2}}\right)^{\alpha+{1\over2}}
\left(\cos{\displaystyle{\theta\over2}}\right)^{\beta+{1\over2}}}\tag4.9\\
E_n^{(\alpha,\beta)}(\cos\theta)&=O(n^{-3/2}),
\endalign
$$
where $\alpha,\beta\in\R$ and $0<\theta<\pi$ are fixed.

We show in Lemma 4.2 below that replacing in $I_{kl}$ the Jacobi polynomials by their 
Darboux approximations leaves the asymptotics of $I_{kl}$ unchanged. This will follow
from the following inequalities due to Szeg\H o \cite{Sz1} (see also 
\cite{Sz2,\, p.197}).

\proclaim{Proposition 4.1 (Szeg\H o \cite{Sz1,(46'),\,(48'),\,{\rm p.77}}\cite{Sz2})} $(a)$ Let
$\alpha\geq-1/2$, $\beta\in\R$ and $\epsilon>0$. Then there exists a constant $A$
depending only on $\alpha$, $\beta$ and $\epsilon$ so that
$$
\left|P_n^{(\alpha,\beta)}(x)\right|\leq\frac{A}
{\sqrt{n}}\frac{1}{(1-x)^{{\alpha\over2}+{1\over4}}}, 
\ \ \ \ \ -1+\epsilon\leq x\leq1.\tag4.10
$$
$(b)$ Let $\alpha,\beta\in\R$ and $\epsilon,c>0$. Then there exists a
constant $B$ depending only on $\alpha$, $\beta$, $\epsilon$ and $c$ so that
$$
\left|E_n^{(\alpha,\beta)}(\cos\theta)\right|\leq\frac{B}{n^{3/2}}
\frac{1}{\theta^{\alpha+{3\over2}}}, 
\ \ \ \ \ cn^{-1}\leq\theta\leq\pi-\epsilon.\tag4.11
$$
\endproclaim

One readily checks that
$$
\cos^{-1}(x)\geq\sqrt{2}\sqrt{1-x}, \ \ \ \ \ 0\leq x\leq1.\tag4.12
$$
Indeed, if $f(x):=\cos^{-1}(x)-\sqrt{2}\sqrt{1-x}$, one has
$f'(x)=(1-x)^{-1/2}(2^{-1/2}-(1+x)^{-1/2})<0$, for $x\in[0,1)$, and as $f(1)=0$, 
(4.12) follows.

Using (4.12) and Proposition 4.1(b) one obtains that for any $\epsilon>0$ there exists
constants $B_\epsilon$ and $N_\epsilon$ depending only on $\alpha$, $\beta$ and 
$\epsilon$ so that
$$
\left|E_n^{(\alpha,\beta)}(x)\right|\leq\frac{B_\epsilon}{n^{3/2}}
\frac{1}{(1-x)^{{\alpha\over2}+{3\over4}}}, 
\ \ \ \ \ {1\over2}\leq x\leq1-\epsilon,\ n\geq N_\epsilon.\tag4.13
$$
By Proposition 4.1(a), there exists a constant $A$ depending only on $\alpha$ and 
$\beta$ so that 
$$
\left|P_n^{(\alpha,\beta)}(x)\right|\leq\frac{A}
{\sqrt{n}}\frac{1}{(1-x)^{{\alpha\over2}+{1\over4}}}, 
\ \ \ \ \ {1\over2}\leq x\leq1.\tag4.14
$$
Using (4.9) and the fact that $\sin(\cos^{-1}(x)/2)=\sqrt{(1-x)/2}$ and 
$\cos(\cos^{-1}(x)/2)=\sqrt{(1+x)/2}$, one obtains that the Darboux approximants
satisfy an inequality of the same form as (4.14): there exists a constant $A'$ such
that
$$
\left|D_n^{(\alpha,\beta)}(x)\right|\leq\frac{A'}
{\sqrt{n}}\frac{1}{(1-x)^{{\alpha\over2}+{1\over4}}}, 
\ \ \ \ \ {1\over2}\leq x\leq1.\tag4.15
$$
By (4.14) and (4.15), 
$E_n^{(\alpha,\beta)}(x)=P_n^{(\alpha,\beta)}(x)-D_n^{(\alpha,\beta)}(x)$
also satisfies an inequality of the same type, so there exists a constant $A''$
depending only on $\alpha$ and $\beta$ so that
$$
\left|E_n^{(\alpha,\beta)}(x)\right|\leq\frac{A''}
{\sqrt{n}}\frac{1}{(1-x)^{{\alpha\over2}+{1\over4}}}, 
\ \ \ \ \ {1\over2}\leq x\leq1.\tag4.16
$$
We are now ready to prove the announced invariance of the asymptotics of
$I_{kl}(R_1,R_2)$ under replacement of the Jacobi polynomials by their Darboux
approximants. Define
$$
J_{kl}(R_1,R_2):=\int_0^1 D_{R_1-k}^{(v_1+k+1,k-v_1-2)}\left(1-{x\over2}\right)
D_{R_2-l}^{(v_2+l,l-v_2-1)}\left(1-{x\over2}\right)x^{k+l+u-1}dx.\tag4.17
$$

\proclaim{Lemma 4.2} Let $r$, $v_1$, $v_2$, $0\leq k\leq v_1+1$ and 
$0\leq l\leq v_2$ be fixed, and
assume that not all of $v_1$, $v_2$, $k$ and $l$ are zero. Then we have
$$
\lim_{R\to\infty}R\left(I_{kl}(R+r,R)-J_{kl}(R+r,R)\right)=0.
$$
\endproclaim

\pf Denoting for simplicity 
$f(x):=P_{R_1-k}^{(v_1+k+1,k-v_1-2)}(x)$, $g(x):=P_{R_2-l}^{(v_2+l,l-v_2-1)}(x)$,
$F(x):=D_{R_1-k}^{(v_1+k+1,k-v_1-2)}(x)$ and $G(x):=D_{R_2-l}^{(v_2+l,l-v_2-1)}(x)$,
we obtain
$$
\align
&\left|I_{kl}(R_1,R_2)-J_{kl}(R_1,R_2)\right|\\
&=\left|\int_0^1\left\{f(1-x/2)g(1-x/2)-F(1-x/2)G(1-x/2)\right\}x^{k+l+u-1}dx\right|\\
&=|\int_0^1\{f(1-x/2)g(1-x/2)-f(1-x/2)G(1-x/2)+f(1-x/2)G(1-x/2)\\ 
&\ \ \ \ \ \ \ \ \ \ \ \ \ \ \ \ \ \ \ \ \ \ \ \ \ \ \ \ \ \ \ \ \ \ \ \ \ \ 
\ \ \ \ \ \ \ \ \ \ \ \ \ \ \ \ \ \ \ \ \ \ \ \ \ -F(1-x/2)G(1-x/2)\}x^{k+l+u-1}dx|
\\
&\leq\int_0^1\left|f(1-x/2)\right|\left|g(1-x/2)-G(1-x/2)\right|x^{k+l+u-1}dx\\
&\ \ \ \ \ \ \ \ \ \ \ \ \ +\int_0^1 \left|G(1-x/2)\right|\left|
f(1-x/2)-F(1-x/2)G(1-x/2)\right|x^{k+l+u-1}dx.
\endalign
$$
Therefore, to prove the Lemma it suffices to show that for $R_1=R+r$ and $R_2=R$,
$$
\lim_{R\to\infty}R\int_0^1\left|P_{R_1-k}^{(v_1+k+1,k-v_1-2)}(1-x/2)\right|
\left|E_{R_2-l}^{(v_2+l,l-v_2-1)}(1-x/2)\right|x^{k+l+u-1}dx=0\tag4.18
$$
and
$$
\lim_{R\to\infty}R\int_0^1\left|E_{R_1-k}^{(v_1+k+1,k-v_1-2)}(1-x/2)\right|
\left|D_{R_2-l}^{(v_2+l,l-v_2-1)}(1-x/2)\right|x^{k+l+u-1}dx=0.\tag4.19
$$
By the Cauchy-Schwarz inequality, for $h$ with integrable square we have
$$
\left(\int_0^1h(x)dx\right)^2\leq\int_0^1h^2(x)dx.\tag4.20
$$
Let $0<\epsilon<1$ be arbitrary. By (4.20) we have
$$
\align
&\left(\int_0^1\left|P_{R_1-k}^{(v_1+k+1,k-v_1-2)}(1-x/2)\right|
\left|E_{R_2-l}^{(v_2+l,l-v_2-1)}(1-x/2)\right|x^{k+l+u-1}dx\right)^2\\
&\leq\int_0^{1-\epsilon}\left|P_{R_1-k}^{(v_1+k+1,k-v_1-2)}(1-x/2)\right|^2
\left|E_{R_2-l}^{(v_2+l,l-v_2-1)}(1-x/2)\right|^2 x^{2(k+l+u-1)}dx\\
&+\int_{1-\epsilon}^{1}\left|P_{R_1-k}^{(v_1+k+1,k-v_1-2)}(1-x/2)\right|^2
\left|E_{R_2-l}^{(v_2+l,l-v_2-1)}(1-x/2)\right|^2 x^{2(k+l+u-1)}dx.\tag4.21
\endalign
$$
Denote the two integrals on the right hand side of (4.21) by $I_1$ and $I_2$,
respectively. By (4.13) and (4.14), we have
$$\align
&I_1\leq\frac{A^2B_\epsilon^2}
{(R_1-k)(R_2-l)^3}\int_0^{1-\epsilon}(x/2)^{-(v_1+k+1+{1\over2})}
(x/2)^{-(v_2+l+{3\over2})}x^{2(k+l+u-1)}dx\\
&\leq\frac{M_\epsilon}{R_1R_2^3}\int_0^{1-\epsilon}x^{v_1+v_2+k+l-1}dx,\tag4.22
\endalign
$$
where $M_\epsilon$ depends on $v_1$, $v_2$, $k$, $l$ and $\epsilon$, but is 
independent of $R_1$ and $R_2$ (here we used that $u=v_1+v_2+2$).
Since by hypothesis $v_1$, $v_2$, $k$ and $l$ are not all equal to zero, the exponent 
of the integrand in (4.22) is nonnegative and we obtain
$$
I_1\leq\frac{M_\epsilon}{R_1R_2^3}.\tag4.23
$$
On the other hand, by (4.14) and (4.16), we have that
$$
\align
&I_2\leq\frac{AA''}{(R_1-k)(R_2-l)}\int_{1-\epsilon}^1
(x/2)^{-(v_1+k+1+{1\over2})}(x/2)^{-(v_2+l+{1\over2})}x^{2(k+l+u-1)}dx\\
&\leq\frac{M'}{R_1R_2}\int_{1-\epsilon}^1x^{v_1+v_2+k+l}dx\\
&\leq\frac{\epsilon M'}{R_1R_2},\tag4.24
\endalign
$$
where $M'$ is independent of $\epsilon$, $R_1$ and $R_2$, depending just on 
$v_1$, $v_2$, $k$ and $l$.

By (4.21), (4.23) and (4.24), for $R_1=R+r$ and $R_2=R$ the first term in (4.21) is
majorized by
$$
I_1+I_2\leq\frac{\epsilon M'}{R^2}+\frac{M_\epsilon}{R^4}\leq\frac{1}{R^2}
\left(\sqrt{\epsilon M'}+\frac{\sqrt{M_\epsilon}}{R}\right)^2.
$$
Extracting the square root we obtain
$$
\align
&\int_0^1\left|P_{R_1-k}^{(v_1+k+1,k-v_1-2)}(1-x/2)\right|
\left|E_{R_2-l}^{(v_2+l,l-v_2-1)}(1-x/2)\right|x^{k+l+u-1}dx\\
&\leq\frac{\sqrt{\epsilon M'}}{R}+\frac{\sqrt{M_\epsilon}}{R^2}.
\endalign
$$
Multiplying the previous inequality by $R$, we obtain that the quantity whose limit is
taken in (4.18) is majorized by $\sqrt{\epsilon M'}+\sqrt{M_\epsilon}/R$. This
quantity can be made arbitrarily small by first choosing $\epsilon$ so as to make 
$\sqrt{\epsilon M'}$ arbitrarily small, and then requiring $R$ to be large enough to
make $\sqrt{M_\epsilon}/R$ arbitrarily small. This proves (4.18).

A similar argument proves (4.19). This completes the proof of the Lemma. \endpf

\mysec{5. The asymptotics of $M_1$}

The following result will be needed several times during the remaining part of the
paper.

\proclaim{Lemma 5.1} 
Let $\alpha(t)$, $k(t)$ and $h(t)$ be complex-valued functions that are real for real
$t$ and analytic in a domain ${\bold T}\subset\C$ containing the interval $(0,1]$.
Let $q>0$ be fixed. Then
$$
\align
\int_0^1t^{Rq}h(t)&\cos[R\alpha(t)+k(t)]dt\\
&=
\frac{h(1)}{R\sqrt{q^2+(\alpha'(1))^2}}\cos\left[R\alpha(1)+k(1)
-\arctan\frac{\alpha'(1)}{q}\right]+O(R^{-2}).\tag5.1
\endalign
$$
\endproclaim

\pf We have
$$
\align
\int_0^1t^{Rq}&h(t)\cos[R\alpha(t)+k(t)]dt\\
&=\frac{1}{2}\int_0^1e^{Rq\ln t}
\left[e^{i(R\alpha(t)+k(t))}+e^{-i(R\alpha(t)+k(t))}\right]h(t)dt\\
&=-\frac{1}{2}\left\{\int_1^0e^{-R[-q\ln t-i\alpha(t)]}e^{ik(t)}h(t)dt
+\int_1^0e^{-R[-q\ln t+i\alpha(t)]}e^{-ik(t)}h(t)dt\right\}.\tag5.2
\endalign
$$
The asymptotics for large $R$ of each of the two integrals on the last line of (5.2) 
can be found by the Laplace method as it is described for example in \S6,
Chapter 4 of \cite{O}.
Indeed, consider the first integral. The only requirement of the hypothesis of 
Theorem 6.1 of \cite{O,\,p.\,125} that
needs to be checked is that the real part of the coefficient of $-R$ in the first 
exponential in the integrand attains its minimum at $t=1$. This indeed holds, since
$-\ln t$ has its minimum at $t=1$.
The relevant quantities are easily found to be, in the notation 
of the quoted theorem of \cite{O}, $\lambda=\mu=1$, $p_0=-q-i\alpha'(1)$,
$q_0=e^{ik(1)}h(1)$ and $p(1)=-i\alpha(1)$. By that theorem, the value of the integral 
is $a_0e^{-Rp(1)}/R+O(R^{-2})$, where $a_0=q_0/(\mu p_0^{\lambda/\mu})$. Therefore, we
obtain by Theorem 6.1 of \cite{O,\,p.\,125} that
$$
\int_1^0e^{-R[-q\ln t-i\alpha(t)]}e^{ik(t)}h(t)dt=
-\frac{e^{i(R\alpha(1)+k(1))}h(1)}{q+i\alpha'(1)}\frac{1}{R}
+O(R^{-2}).\tag5.3
$$
The two integrals on the last line of (5.2) are complex conjugates, so the 
coefficients of their asymptotic expansions are also complex conjugates. We obtain
from (5.3) that
$$
\int_1^0e^{-R[-q\ln t+i\alpha(t)]}e^{-ik(t)}h(t)dt=
-\frac{e^{-i(R\alpha(1)+k(1))}h(1)}{q-i\alpha'(1)}\frac{1}{R}
+O(R^{-2}).\tag5.4
$$
However, it is readily checked that
$$
\frac{e^{i\varphi}}{q+ib}+\frac{e^{-i\varphi}}{q-ib}=\frac{2}{\sqrt{q^2+b^2}}
\cos(\varphi-\theta),\tag5.5
$$
where $\theta=\arctan(b/q)$. 
By (5.2)--(5.5) we obtain the statement of the Lemma. \endpf

\proclaim{Lemma 5.2} For fixed $r$, $v_1$ and $v_2$ we have
$$
\align
&I_{v_1+1,v_2}(R+r,R)\\
&\ \ \ 
=\frac{(-4)^{v_1+v_2+1}}{2R\pi}\int_0^1\left(\frac{4-x}{x}\right)^{1/2}x^{v_1+v_2+1}
\cos\left[r\cos^{-1}\left(1-{x\over2}\right)\right]dx
+o(R^{-1}).\tag5.6
\endalign
$$
\endproclaim

\pf By Lemma 4.2, it is enough to show that $J_{v_1+1,v_2}(R+r,R)$ has the asymptotics
given by (5.6). Using (4.17), (4.9) and the fact that 
$\sin(\cos^{-1}(1-x/2)/2)=(x/4)^{1/2}$,
$\cos(\cos^{-1}(1-x/2)/2)=((4-x)/4)^{1/2}$ and $\cos(z-n\pi)=(-1)^n\cos(z)$, 
one obtains that
$$
\align
&J_{v_1+1,v_2}(R+r,R)\\
&=\frac{(-4)^{v_1+v_2+1}}{\pi\sqrt{R+r-v_1-1}\sqrt{R-v_2}}
\int_0^1\left(\frac{4-x}{x}\right)^{1/2}x^{v_1+v_2+1}
\cos\left[(R+r)\cos^{-1}\left(1-{x\over2}\right)-{\pi\over4}\right]\\
&\ \ \ \ \ \ \ \ \ \ \ \ \ \ \ \ \ \ \ \ \ \ \ \ \ \ \ \ \ \ \ \ \ \ \ \ \ \ \ 
\times\cos\left[R\cos^{-1}\left(1-{x\over2}\right)-{\pi\over4}\right]dx.\tag5.7
\endalign
$$
Converting the product of cosines into a sum, we obtain that
$$
J_{v_1+1,v_2}(R+r,R)=\frac{(-4)^{v_1+v_2+1}}{2\pi\sqrt{R+r-v_1-1}\sqrt{R-v_2}}
(J_1+J_2),\tag5.8
$$
where
$$
\align
J_1&=\int_0^1\left(\frac{4-x}{x}\right)^{1/2}x^{v_1+v_2+1}
\cos\left[r\cos^{-1}\left(1-{x\over2}\right)\right]dx,\tag5.9\\
J_2&=\int_0^1\left(\frac{4-x}{x}\right)^{1/2}x^{v_1+v_2+1}
\cos\left[(2R+r)\cos^{-1}\left(1-{x\over2}\right)-{\pi\over2}\right]dx.
\endalign
$$
By Lemma 5.1, $J_2=O(1/R)$. Therefore, (5.7)--(5.9) imply that
the asymptotics of $J_{v_1+1,v_2}(R+r,R)$ is given by the right hand side of (5.6). As
noted in the beginning of the proof, this completes the proof of the Lemma. \endpf

We are now ready to give the asymptotics of $M_1$.

\proclaim{Proposition 5.3} 
For fixed $r$, $v_1$ and $v_2$, we have
$$
M_1(R+r,R)=\frac{1}{R^3}\frac{1}{\pi}\int_0^1\left(\frac{4-x}{x}\right)^{1/2}x^{u-1}
\cos\left[r\cos^{-1}\left(1-{x\over2}\right)\right]dx
+o(R^{-3}),\tag5.10
$$
where $u=v_1+v_2+2$.
\endproclaim

\pf Using the bounds (4.14) for the Jacobi polynomials, we obtain from (4.8) that
$$
\align
&I_{kl}(R_1,R_2)\leq\frac{M}{\sqrt{R_1-k}\sqrt{R_2-l}}
\int_0^1 \frac{1}{\left({\displaystyle {x\over2}}\right)^{{v_1+k+1\over2}+{1\over4}}}
\frac{1}{\left({\displaystyle {x\over2}}\right)^{{v_2+l\over2}+{1\over4}}}
x^{v_1+v_2+k+l+1}dx\\
&\leq\frac{M'}{\sqrt{R_1}\sqrt{R_2}}\int_0^1 x^{v_1+v_2+k+l\over2}dx,\tag5.11
\endalign
$$
where the constants $M$ and $M'$ depend just on $v_1$, $v_2$, $k$ and $l$. 
Therefore, we obtain that
$$
I_{kl}(R+r,R)=O\left({1\over R}\right),\ \ \ \ \ 0\leq k\leq v_1+1, 0\leq l\leq v_2
.\tag 5.12
$$
Consider now the representation of $M_1$ given by (4.7). For $R_1=R+r$ and $R_2=R$, the 
coefficient of each $I_{kl}(R_1,R_2)$ with $(k,l)\neq(v_1+1,v_2)$ is $O(R^{-3})$. Thus,
by (5.12) we obtain from (4.7) that
$$
M_1(R+r,R)=\frac{c_{v_1+1,v_2}}{R^2}I_{v_1+1,v_2}(R+r,R)+O\left({1\over R^4}\right).
\tag5.13
$$
By (3.8), (4.2), (4.3), (4.5) and (4.6) we obtain after simplifications that
$$
c_{v_1+1,v_2}=\frac{2}{(-4)^{v_1+v_2+1}}.
$$
Substituting this value into (5.13) and using Lemma 5.2 we obtain the statement of the
Proposition. \endpf

\mysec{6. The asymptotics of $M_a$, $M_c$ and $M_{ac}$}

Our analysis of the asymptotics of $M_1$ can be repeated for the remaining double sums
$M_a$, $M_c$ and $M_{ac}$. We obtain the following result.

\proclaim{Proposition 6.1} 
For fixed $r$, $v_1$ and $v_2$, we have
$$
\align
M_a(R+r,R)&=-\frac{1}{R^2}\frac{1}{\pi}\int_0^1x^{u-1}
\cos\left[r\cos^{-1}\left(1-{x\over2}\right)-{\pi\over2}\right]dx
+o(R^{-2})\tag6.1\\ \\
M_c(R+r,R)&=\frac{1}{R^2}\frac{1}{\pi}\int_0^1x^{u-1}
\cos\left[r\cos^{-1}\left(1-{x\over2}\right)-{\pi\over2}\right]dx
+o(R^{-2})\tag6.2\\ \\
M_{ac}(R+r,R)&=\frac{1}{R}\frac{1}{\pi}\int_0^1\left(\frac{4-x}{x}\right)^{-1/2}x^{u-1}
\cos\left[r\cos^{-1}\left(1-{x\over2}\right)\right]dx
+o(R^{-1})\tag6.3
\endalign
$$
where $u=v_1+v_2+2$.
\endproclaim

\pf By (4.1) we can express the $_3F_2$'s of (3.9) and (3.10) as
$$
\align
\!\!\!\!\!\!\!\!\!
&_3F_2\!\left[\matrix{-R_1+1,R_1+1,v_1+{5\over2}}\\{{3\over2},v_1+3}\endmatrix;
{\displaystyle {x\over4}}\right]=\frac{1}{(-v_1-{3\over2})_{v_1+1}}
\sum_{k=0}^{v_1+1} (-1)^k{v_1+1\choose k}(-v_1-{3\over2})_{v_1+1-k}\\
&\ \ \ \ \ \ \ \ \ \ \ \ \ \ \ \ \ \ \ \ \ \ \ \ \ \ \ \
\times\frac{(-R_1+1)_k(R_1+1)_k}{(v_1+3)_k}\frac{x^k}{4^k}\, 
_2F_1\!\left[\matrix{-R_1+1+k,R_1+1+k}\\{v_1+3+k}\endmatrix;{\displaystyle {x\over4}}
\right]\tag6.4
\endalign
$$
and
$$
\align
&_3F_2\!\left[\matrix{-R_2+1,R_2+1,v_2+{5\over2}}\\{{5\over2},v_2+2}\endmatrix;
{\displaystyle {x\over4}}\right]=\frac{1}{(-v_2-{3\over2})_{v_2}}
\sum_{l=0}^{v_2} (-1)^l{v_2\choose l}(-v_2-{3\over2})_{v_2-l}\\
&\ \ \ \ \ \ \ \ \ \ \ \ \ \ \ \ \ \ \ \ \ \ \ \ \ \  
\times\frac{(-R_2+1)_l(R_2+1)_l}{(v_2+2)_l}\frac{x^l}{4^l}\, 
_2F_1\!\left[\matrix{-R_2+1+l,R_2+1+l}\\{v_2+2+l}\endmatrix;{\displaystyle {x\over4}}
\right].\tag6.5
\endalign
$$
By (4.4), the resulting $_2F_1$'s are expressed in terms of Jacobi polynomials as
$$
_2F_1\!\left[\matrix{-R_1+1+k,R_1+1+k}\\{v_1+3+k}\endmatrix;
{\displaystyle {x\over4}}\right]=\frac{(R_1-k-1)!\,(v_1+k+2)!}{(R_1+v_1+1)!}
P_{R_1-k-1}^{(v_1+k+2,k-v_1-1)}\!\left(1-{x\over2}\right)\tag6.6
$$
and
$$
_2F_1\!\left[\matrix{-R_2+1+l,R_2+1+l}\\{v_2+2+l}\endmatrix;
{\displaystyle {x\over4}}\right]=\frac{(R_2-l-1)!\,(v_2+l+1)!}{(R_2+v_2)!}
P_{R_2-l-1}^{(v_2+l+1,l-v_2)}\left(1-{x\over2}\right).\tag6.7
$$
Substituting the expansions (6.4) and (4.3) and the formulas (6.6) and (4.6) into the
integral representation (3.11) of $M_a$, we obtain that
$$
M_a=\sum_{k=0}^{v_1+1}\sum_{l=0}^{v_2} \frac{c_{kl}'R_1}{R_2}
\frac{(-R_1+1)_k(R_1+1)_k}{(R_1-k)_{v_1+k+2}}
\frac{(-R_2)_l(R_2)_l}{(R_2-l+1)_{v_2+l}} I_{kl}'(R_1,R_2),\tag6.8
$$
where $c_{kl}'$ depends only on $k$, $l$, $v_1$ and $v_2$, and for $0\leq k\leq v_1+1$ 
and $0\leq l\leq v_2$,
$$
I_{kl}'(R_1,R_2):=\int_0^1 P_{R_1-k-1}^{(v_1+k+2,k-v_1-1)}\left(1-{x\over2}\right)
P_{R_2-l}^{(v_2+l,l-v_2-1)}\left(1-{x\over2}\right)x^{k+l+u}dx.\tag6.9
$$
Similarly, we get from (3.12) that
$$
M_c=\sum_{k=0}^{v_1+1}\sum_{l=0}^{v_2} \frac{c_{kl}''R_2}{R_1}
\frac{(-R_1)_k(R_1)_k}{(R_1-k+1)_{v_1+k+1}}
\frac{(-R_2+1)_l(R_2+1)_l}{(R_2-l)_{v_2+l+1}} I_{kl}''(R_1,R_2),\tag6.10
$$
where $c_{kl}''$ depends only on $k$, $l$, $v_1$ and $v_2$, and for $0\leq k\leq v_1+1$ 
and $0\leq l\leq v_2$,
$$
I_{kl}''(R_1,R_2):=\int_0^1 P_{R_1-k}^{(v_1+k+1,k-v_1-2)}\left(1-{x\over2}\right)
P_{R_2-l-1}^{(v_2+l+1,l-v_2)}\left(1-{x\over2}\right)x^{k+l+u}dx.\tag6.11
$$
An analogous calculation yields from (3.13) that 
$$
M_{ac}=\sum_{k=0}^{v_1+1}\sum_{l=0}^{v_2} c_{kl}'''R_2R_1
\frac{(-R_1+1)_k(R_1+1)_k}{(R_1-k)_{v_1+k+2}}
\frac{(-R_2+1)_l(R_2+1)_l}{(R_2-l)_{v_2+l+1}} I_{kl}'''(R_1,R_2),\tag6.12
$$
where $c_{kl}'''$ depends only on $k$, $l$, $v_1$ and $v_2$, and for 
$0\leq k\leq v_1+1$ and $0\leq l\leq v_2$,
$$
I_{kl}'''(R_1,R_2):=\int_0^1 P_{R_1-k-1}^{(v_1+k+2,k-v_1-1)}\left(1-{x\over2}\right)
P_{R_2-l-1}^{(v_2+l+1,l-v_2)}\left(1-{x\over2}\right)x^{k+l+u+1}dx.\tag6.13
$$
As seen in the proof of Lemma 5.2 for the case of $I_{kl}$, the bounds (4.14) imply that
for fixed $r$, $v_1$, $v_2$ and fixed $0\leq k\leq v_1+1$ and $0\leq l\leq v_2$, the
integrals $I_{kl}'(R+r,R)$, $I_{kl}''(R+r,R)$ and $I_{kl}'''(R+r,R)$ are $O(1/R)$.
Indeed, the only change from that case is that now the parameters of the 
Jacobi polynomials $P_n^{(\alpha,\beta)}(x)$ that occur are slightly changed. The key
fact needed to prove (5.12) was that the exponent of $x$ in the last integral of
(5.11) was nonnegative. However, the analogous exponents for the case of
$I_{k,l}'$, $I_{k,l}''$ and $I_{k,l}'''$ are readily seen to be nonnegative 
as well (by (4.14), this exponent goes down half a unit for each unit of increase 
in the $\alpha$-parameter of the Jacobi polynomials that occur; the $\alpha$-parameters
of the pairs
of Jacobi polynomials $P_n^{(\alpha,\beta)}(x)$ appearing in $I_{kl}'$, $I_{kl}''$ and
$I_{kl}'''$  are increased by $(1,0)$, $(0,1)$ and $(1,1)$, respectively; the increase
in the exponent of $x$ in (3.11)--(3.13), namely $1$, $1$, and $2$ units, respectively,
makes up for the decrease due to the change in the $\alpha$-parameters).

Using this, it can be shown that, just as it was the case for the expansion (4.7) of 
$M_1$, the asymptotics of the double
sums (6.8), (6.10) and (6.12) for $r$, $v_1$, $v_2$ fixed and $R_1=R+r$, $R_2=R$, 
$R\to\infty$ are given by the contribution of the terms with $(k,l)=(v_1+1,v_2)$. 

To work these out, note first that analogs of Lemma 4.2 hold for $I_{v_1+1,v_2}'$,
$I_{v_1+1,v_2}''$ and $I_{v_1+1,v_2}'''$, with $J_{v_1+1,v_2}'$, $J_{v_1+1,v_2}''$ and 
$J_{v_1+1,v_2}'''$ defined by replacing the Jacobi polynomials in the integrands of
the $I$-integrals by their Darboux approximants. Indeed, the only difference from the
calculations in the proof of Lemma 4.2 is that now the $\alpha$-parameters of the pairs
of Jacobi polynomials $P_n^{(\alpha,\beta)}(x)$ appearing in the $I$-integrals are 
increased by $(1,0)$, $(0,1)$ and $(1,1)$, respectively. However, just as was the case 
in the previous paragraph, the increase in the exponent of $x$ in (6.9), (6.11) and 
(6.13) (of $1$, $1$ and $2$ units, respectively) compensates the decrease due to the
change in the $\alpha$-parameters. Therefore, arguments parallel to the ones in the
proof of Lemma 4.2 lead to analogs of the majorizations (4.22) and (4.24) that
maintain the key feature of having non-negative exponents of $x$ in the integrand, and
thus prove the claimed analogs of Lemma 4.2.

Second, using these analogs of Lemma 4.2, one can easily deduce analogs of Lemma 5.2,
yielding
$$
\align
&I'_{v_1+1,v_2}(R+r,R)\\
&\ \ \ 
=-\frac{(-4)^{v_1+v_2+2}}{2R\pi}\int_0^1x^{v_1+v_2+1}
\cos\left[(r-1)\cos^{-1}\left(1-{x\over2}\right)\right]dx
+o(R^{-1})\tag6.14\\ \\
&I''_{v_1+1,v_2}(R+r,R)\\
&\ \ \ 
=\frac{(-4)^{v_1+v_2+2}}{2R\pi}\int_0^1x^{v_1+v_2+1}
\cos\left[(r+1)\cos^{-1}\left(1-{x\over2}\right)\right]dx
+o(R^{-1})\tag6.15\\ \\
&I'''_{v_1+1,v_2}(R+r,R)\\
&\ \ \ 
=\frac{(-4)^{v_1+v_2+3}}{2R\pi}\int_0^1\left(\frac{4-x}{x}\right)^{-1/2}x^{v_1+v_2+1}
\cos\left[r\cos^{-1}\left(1-{x\over2}\right)\right]dx
+o(R^{-1}).\tag6.16
\endalign
$$
And third, replacing in (3.11)--(3.13) the expansions (4.2), (4.3), (6.4) and  (6.5)
and formulas (4.5), (4.6), (6.6) and (6.7), the constants $c_{v_1+1,v_2}'$,
$c_{v_1+1,v_2}''$ and $c_{v_1+1,v_2}'''$ are found, after simplifications, to be 
$$
\align
&c_{v_1+1,v_2}'=\frac{2}{(-4)^{v_1+v_2+2}}\\ 
&c_{v_1+1,v_2}''=-\frac{2}{(-4)^{v_1+v_2+2}}\\ 
&c_{v_1+1,v_2}'''=\frac{2}{(-4)^{v_1+v_2+3}}.
\endalign
$$
Substituting these and (6.14)--(6.16) into (6.8), (6.10) and (6.12) we obtain the
statements (6.1)--(6.3) of the Proposition. \endpf

\mysec{7. The asymptotics of the correlation $\omega(r,u)$}

Substituting the asymptotics of the double sums $M_1$, $M_a$, $M_c$ and $M_{ac}$ given
by Propositions 5.3 and 6.1 into the formula (3.3), we obtain the following result.

\proclaim{Proposition 7.1} For fixed $r$, $v_1$ and $v_2$, we have
$$
\omega_b(R+r,v_1;R,v_2)=\frac{1}{4\pi^2}|S_1S_{ac}+S_aS_c|+o(R^{-1}),
$$
where
$$
\align
S_1&=\int_0^1\left(\frac{4-x}{x}\right)^{1/2}x^{u-1}
\cos\left[r\cos^{-1}\left(1-{x\over2}\right)\right]dx\\
S_a&=\int_0^1x^{u-1}
\cos\left[r\cos^{-1}\left(1-{x\over2}\right)-{\pi\over2}\right]dx\\
S_c&=\int_0^1x^{u-1}
\cos\left[r\cos^{-1}\left(1-{x\over2}\right)-{\pi\over2}\right]dx\\
S_{ac}&=\int_0^1\left(\frac{4-x}{x}\right)^{-1/2}x^{u-1}
\cos\left[r\cos^{-1}\left(1-{x\over2}\right)\right]dx
\endalign
$$
and $u=v_1+v_2+2$.
\endproclaim
\flushpar
{\smc Remark 7.2.} By the above result, for fixed $r$, $v_1$ and $v_2$ the 
asymptotics of $\omega_b(R+r,v_1;R,v_2)$ as $R\to\infty$ depends only on the sum
$v_1+v_2$, and not individually on $v_1$ and $v_2$. This is consistent with the
expectation that the quadromer correlation at the center should depend only on the
separation vector $(r,u)$.


\medskip
In the statement of Theorem 1.1, the coordinates of the separation vector $(r,u)$ are
related by $u=qr+c$, where $q\geq0$ and $c$ are fixed rational numbers. When $q>0$,
the asymptotics of $\omega_b(R+r,v_1;R,v_2)$ as $R\to\infty$ can be obtained from
Lemma 5.1. To handle the case $q=0$ we need the following result.

\proclaim{Lemma 7.3} 
Let $\alpha$, $k$ and $h$ be real-valued functions that are analytic in an open
interval containing $(0,1]$. Assume $\alpha'(t)>0$ in $(0,1)$ and 
$\lim_{t\to0^+}h(t)/\alpha'(t)=0$.
Then
$$
\align
\int_0^1h(t)&\cos[R\alpha(t)+k(t)]dt\\
&=
\frac{h(1)}{R\alpha'(1)}\cos\left[R\alpha(1)+k(1)
-{\pi\over2}\right]+O(R^{-2}).\tag7.1
\endalign
$$
\endproclaim

\pf As in the proof of Lemma 5.1, express the integrand in terms of exponentials as
$$
\align
\int_0^1&h(t)\cos[R\alpha(t)+k(t)]dt\\
&=\frac{1}{2}\left\{\int_0^1e^{iR\alpha(t)}e^{ik(t)}h(t)dt
+\int_0^1e^{-iR\alpha(t)}e^{-ik(t)}h(t)dt\right\}.\tag7.2
\endalign
$$
Consider the first integral on the right hand side of (7.2). Make the change
of variables $y=1-t$ to obtain
$$
\int_0^1e^{iR\alpha(t)}e^{ik(t)}h(t)dt=
\int_0^1e^{-iR\gamma(t)}\delta(t)dt,\tag7.3
$$
where $\gamma(y)=-\alpha(1-y)$ and $\delta(y)=e^{ik(1-y)}h(1-y)$. 
These functions $\gamma(y)$ and $\delta(y)$ are readily checked to satisfy the
conditions in the hypothesis of Theorem 13.2 of \cite{O,\,p.\,102} with $i$
replaced by $-i$ throughout (clearly, by complex conjugation, the statement of the 
quoted theorem remains true when $i$ is replaced by $-i$ throughout; we need to apply
this modified version of the quoted theorem because its hypothesis requires  
$\gamma(y)'>0$; compare with the beginning of \S13.1 of \cite{O}).
Since the functions $\alpha(t)$, $k(t)$ and $h(t)$ of (7.1) are analytic at $t=1$, it
follows that $\gamma(y)$ and $\delta(y)$ are analytic at $y=0$. Therefore, the
exponents $\lambda$ and $\mu$ of (13.02) \cite{O} are both equal to 1. Thus,
Theorem 13.2 of \cite{O} is applicable and, since we are assuming
$\lim_{t\to0^+}h(t)/\alpha'(t)=0$, it yields
$$
\int_0^1e^{-iR\gamma(t)}\delta(t)dt=\frac{\delta(0)}{\gamma'(0)}
\frac{e^{-iR\gamma(0)}}{iR}+o(R^{-1})\tag7.4
$$
(since $i$ is now replaced by $-i$ throughout Theorem 13.2 of \cite{O}). By
(7.3), (7.4) and the definition of $\gamma(t)$ and $\delta(t)$ we obtain that
$$
\int_0^1e^{iR\alpha(t)}e^{ik(t)}h(t)dt=
\frac{e^{ik(1)}h(1)}{\alpha'(1)}\frac{e^{iR\alpha(1)}}{iR}+o(R^{-1})\tag7.5
$$
The two integrals on the right hand side of (7.2) are complex conjugates, so the 
coefficients of their asymptotic expansions are also complex conjugates. We obtain
from (7.5) that
$$
\int_0^1e^{-iR\alpha(t)}e^{-ik(t)}h(t)dt=
-\frac{e^{-ik(1)}h(1)}{\alpha'(1)}\frac{e^{-iR\alpha(1)}}{iR}+o(R^{-1})\tag7.6
$$
By (7.2)--(7.6) we obtain the statement of the Lemma. \endpf

We are now ready to prove our main result.

{\it Proof of Theorem 1.1.} Let $u=qr+c$, where $q\geq0$ and $q,c\in\Q$ are fixed.
Then the integral $S_1$ of Proposition 7.1 becomes
$$
S_1=\int_0^1\left(\frac{4-t}{t}\right)^{1/2}t^{qr+c-1}
\cos\left[r\cos^{-1}\left(1-{t\over2}\right)\right]dt.
$$
For $q>0$, we can apply Lemma 5.1 with $h(t)=t^{c-3/2}(4-t)^{1/2}$, 
$\alpha(t)=\cos^{-1}(1-t/2)$ and $k(t)=0$ to obtain 
$$
S_1=\frac{1}{r}\frac{\sqrt{3}}{\sqrt{q^2+1/3}}
\cos\left({r\pi\over3}-\arctan{1\over q\sqrt{3}}\right)+o(r^{-1}),\ \ \ \ \ q>0.\tag7.7
$$
If $q=0$, we have $u=c$ and $S_1$ becomes
$$
S_1=\int_0^1\left(\frac{4-t}{t}\right)^{1/2}t^{c-1}
\cos\left[r\cos^{-1}\left(1-{t\over2}\right)\right]dt.
$$
This has the form of the integral in (7.1), with $h(t)=t^{c-1/2}(4-t)^{1/2}$, 
$\alpha(t)=\cos^{-1}(1-t/2)$ and $k(t)=0$. Clearly, $\alpha'(t)=(4t-t^2)^{-1/2}>0$
in $(0,1)$, and it is readily checked that $\lim_{t\to0^+}h(t)/\alpha'(t)=0$ 
(in checking this we need to use the fact that $c=u\geq1$). Therefore, Lemma 7.2 yields
$$
S_1=\frac{3}{r}
\cos\left({r\pi\over3}-{\pi\over2}\right)+o(r^{-1}),\ \ \ \ \ q=0.\tag7.8
$$
A conceptual way of viewing (7.7) and (7.8) together is to say that (7.7) also holds
in the limit $q\to0^+$.
Similar applications of Lemma 5.1 to the integrals $S_a$, $S_c$ and
$S_{ac}$ of Proposition 7.1 yield, for $q>0$, that 
$$
\align
S_a&=\frac{1}{r}\frac{1}{\sqrt{q^2+1/3}}
\cos\left({r\pi\over3}-{\pi\over2}-\arctan{1\over q\sqrt{3}}\right)+o(r^{-1})\tag7.9\\
S_c&=\frac{1}{r}\frac{1}{\sqrt{q^2+1/3}}
\cos\left({r\pi\over3}-{\pi\over2}-\arctan{1\over q\sqrt{3}}\right)+o(r^{-1})\tag7.10\\
S_{ac}&=\frac{1}{r}\frac{1}{\sqrt{3}\sqrt{q^2+1/3}}
\cos\left({r\pi\over3}-\arctan{1\over q\sqrt{3}}\right)+o(r^{-1}),\tag7.11
\endalign
$$
the formulas also holding, by Lemma 7.3, in the limit $q\to0^+$.

By (1.4), Proposition 7.1, (7.7), (7.9)--(7.11) and the fact that the latter four
relations hold also in the limit $q\to0^+$, we obtain
$$
\align
\omega(u&,r)=\frac{1}{4\pi^2}|S_1S_{ac}+S_aS_c|\\
&=\left|\frac{1}{4\pi^2r^2(q^2+1/3)}
\left\{\cos^2\!\!\left({r\pi\over3}-\arctan{1\over q\sqrt{3}}\right)
+\sin^2\!\!\left({r\pi\over3}-\arctan{1\over q\sqrt{3}}\right)\right\}+o(r^{-2})\right|\\
&=\left|\frac{3}{4\pi^2(3q^2r^2+r^2)}+o(r^{-2})\right|\\
&=\frac{3}{4\pi^2(r^2+3u^2)}+o(r^{-2}).
\endalign
$$
This proves Theorem 1.1. \endpf

\medskip
\flushpar
{\bf Acknowledgments.} I would like to thank Jeff Geronimo for useful discussions and for pointing 
out to me Laplace's method for the asymptotics of integrals, and the referee for the careful reading
of the manuscript and helpful suggestions.
 
\mysec{References}
{\openup 1\jot \frenchspacing\raggedbottom
\roster
\myref{C} 
  M. Ciucu, Plane partitions I: A generalization of MacMahon's
formula, preprint (available at the Los Alamos 
archive, at http://arxiv.org/ps/math.CO/9808017).

\myref{FS} 
  M. E. Fisher and J. Stephenson, Statistical mechanics of dimers on a plane 
lattice. II. Dimer correlations and monomers, {\it Phys. Rev. (2)} {\bf 132} (1963),
1411--1431.

\myref{GR}
  G. Gasper and M. Rahman, ``Basic hypergeometric series,'' Cambridge University Press, Cambridge, 
1990.

\myref{GV}
  I. M. Gessel and X. Viennot, Binomial determinants, paths, and hook length formulae, {\it Adv.
in Math.} {\bf 58} (1985), 300--321.


\myref{H} 
  R. E. Hartwig, Monomer pair correlations, {\it J. Mathematical Phys.} {\bf 7}
(1966), 286--299.

\myref{O}
  F. W. J. Olver, Asymptotics and special functions, Academic Press, New York, 1974.


\myref{Sz1} G. Szeg\H o, Asymptotische Entwicklungen der Jacobischen Polynome, {\it Schriften 
K\"o\-nigs\-ber\-ger Gel. Ges.}, (1933), 35--112; re\-prin\-ted in R. Askey (ed.), 
``Collected papers of G\'abor Szeg\H o,'' vol. 2, 399--478, Birkh\"auser, Boston, 1982.

\myref{Sz2}
  G. Szeg\H o, ``Orthogonal polynomials,'' American Mathematical Society,  
Providence, R.I., 1975.

\endroster\par}

\enddocument